\documentclass[final,seceqn,secthm]{elsart}
\usepackage{amsmath,amsbsy,amssymb,amscd}
\usepackage{amsfonts}
\usepackage{fancyhdr}
\usepackage{mathtools}
\usepackage{latexsym}
\usepackage{color}
\usepackage[mathscr]{euscript}
\numberwithin{equation}{section}
\DeclareMathOperator{\diag}{diag}		
\DeclareMathOperator{\sgn}{sgn}			
\DeclareMathOperator{\sinc}{sinc}		
\DeclareMathOperator{\spann}{span}		
\DeclareMathOperator{\supp}{supp}		
\newcommand{\C}{{\mathbb C}}
\newcommand{\N}{{\mathbb N}}
\newcommand{\R}{{\mathbb R}}

\newcommand{\eK}{{\mathscr K}}
\newcommand{\fA}{{\mathfrak A}}
\renewcommand{\d}{\,{\mathrm d}}
\newcommand{\eps}{\varepsilon}
\newcommand{\tm}{\times}
\renewcommand{\leq}{\leqslant}
\newcommand{\Cd}{C(\Omega)^d}
\renewcommand{\le}{\preceq}
\newcommand{\lle}{\prec}
\newcommand{\lld}{\prec\!\!\!\prec}
\newcommand{\intoo}[1]{\left(#1\right)}				
\newcommand{\set}[1]{\left\{#1\right\}}				
\newcommand{\abs}[1]{\left|#1\right|}				
\newcommand{\norm} [1]{\left\|#1\right\|}			
\newcommand{\iprod}[1]{\left\langle#1\right\rangle}
\newcommand{\sprod}[1]{\langle\!\!\langle#1\rangle\!\!\rangle}		
\newcommand{\fall}{\quad\text{for all }}
\newcommand{\cref}[1]{Cor.~\ref{#1}}

\newcommand{\eref}[1]{Ex.~\ref{#1}}
\newcommand{\fref}[1]{Fig.~\ref{#1}}
\newcommand{\lref}[1]{Lemma~\ref{#1}}
\newcommand{\rref}[1]{Rem.~\ref{#1}}
\newcommand{\tref}[1]{Thm.~\ref{#1}}
\allowdisplaybreaks
\begin{document}
\begin{frontmatter}
\title{Positivity and discretization of\\ Fredholm integral operators}
\author[mnr]{Magdalena Nockowska-Rosiak}
\ead{magdalena.nockowska@p.lodz.pl}
\address[mnr]{Institute of Mathematics, Lodz University of Technology,\\
al.\ Politechniki 8, 93-590 \L{}{\'o}d{\'z}, Poland}
\author[cp]{Christian P\"otzsche}
\ead{christian.poetzsche@aau.at}
\address[cp]{Institut f\"ur Mathematik, Universit\"at Klagenfurt,\\
Universit\"atsstra{\ss}e~65--67, 9020 Klagenfurt, Austria}
\begin{abstract}
	We provide sufficient conditions for vector-valued Fredholm integral operators and their commonly used spatial discretizations to be positive in terms of an order relation induced by a corresponding order cone. It turns out that reasonable Nystr\"om methods preserve positivity. Among the projection methods, persistence is obtained for the simplest ones based on polynomial, piecewise linear or specific cubic interpolation (collocation), as well as for piecewise constant basis functions in a Bubnov-Galerkin approach. However, for semi-discreti\-zations using quadratic splines or $\sinc$-collocation we demonstrate that positivity is violated. Our results are illustrated in terms of eigenpairs for Krein-Rutman operators and form the basis of corresponding investigations for nonlinear integral operators. 
\end{abstract}
\end{frontmatter}
\maketitle
\section{Introduction}
Fredholm integral operators canonically arise in numerous applications ranging from a classical fixed-point formulation of linear elliptic boundary value problems \cite{amann:76} to Fr{\'e}chet derivatives of recent models in theoretical ecology \cite[pp.~23ff]{lutscher:19}. Often their Green's function resp.\ kernel is positive. This is understood that such a matrix-valued function preserves an order relation induced by an order cone. For instance, dispersal kernels used in theoretical ecology preserve an order relation in order to capture predator-prey or symbiotic relationships between different species, or differential operators satisfy a maximum principle which transfers to positivity of their inverses via the Green's function. Moreover, for example \cite{erbe:hu:96} provides sufficient conditions that Hammerstein operators can be transformed into an order-preserving form. 

Dealing with positive operators has several merits. First, positivity allows an application of the Krein-Rutman theorem \cite[p.~226, Thm.~19.2]{deimling:85}, \cite[p.~290, Thm.~7.C]{zeidler:93} or \cite{li:jia:21} with profound consequences to the bifurcation behavior of related nonlinear problems. Second, \cite{anselone:lee:74} gives results on the distribution of secondary eigenvalues. Finally, criteria for various types of positivity in linear operators are fundamental to obtain corresponding results for nonlinear integral operators of Urysohn or Hammerstein type \cite{nockowska:poetzsche:1,nockowska:poetzsche:2}. Beyond providing sufficient conditions for positivity of Fredholm integral operators, in this paper we also discuss the question whether this property is preserved under discretizations. The latter aspect is crucial in numerical simulations and computations in order to preserve or capture properties of the original problem. 

More detailled, the texture of this paper is as follows: After introducing our notation, Sect.~\ref{sec2} provides sufficient conditions on kernels such that the associated Fredholm operators on the continuous or integrable functions are positive over domains of finite measure. In essence these properties carry over from the kernels having values in $\R^d$ to the integral operators mapping into a space of $\R^d$-valued functions. Persistence issues under full spatial discretization of Nystr\"om type are tackled in Sect.~\ref{sec3}. It suffices to assume positive quadrature weights for the sake of positivity. This property is satisfied for a large class of integration methods \cite{davis:rabinowitz:84,hammerlin:hoffmann:91} and, besides guaranteeing well-posedness and computational stability \cite{hackbusch:95,hackbusch:14}, provides another reason for the use of positive weights in numerical quadrature. This convenient prospect changes in Sect.~\ref{sec4} when addressing semi-discretizations of projection type. Here the situation is more subtle and the class of positivity-preserving schemes appears to be rather small. Among the collocation methods we indeed observe that polynomial and piecewise linear collocation preserves positivity, while popular other methods (e.g.\ quadratic splines, $\sinc$ methods) do not. Among the Bubnov-Galerkin methods, at least piecewise constant approximation works. Nevertheless, based on positivity-preserving interpolation methods \cite{sakai:schmidt:89,schmidt:hess:88,zhu:18} it appears to be possible to construct projection operators, which are positive at least on subsets of the state space. Pointing towards applications, Sect.~\ref{sec5} discusses order-preserving properties of typical dispersal kernels and illustrates our results by means of a failure in the Krein-Rutman theorem when using quadrature schemes having negative weights. Moreover, we compare numerical errors for various projection methods when approximating the dominant eigenvalue of a positive operator. For the reader's convenience, an appendix collects basic results on cones in Banach spaces and positive operators, an explicit formula for the inverse of tridiagonal matrices and lists the quadrature rules in Nystr\"om methods throughout the text. 

\textbf{Notation}: 
We abbreviate $\R_+:=[0,\infty)$ for the nonnegative reals, $\delta_{ij}\in\set{0,1}$ is the Kronecker symbol and the Euclidean inner product in $\R^d$ is given by $\sprod{x,y}:=\sum_{j=1}^dx_jy_j$ for $x,y\in\R^d$. 

Norms on finite-dimensional spaces are denoted by $\abs{\cdot}$ and unless otherwise stated we use the Euclidean norm. Given a matrix $T\in\R^{d\tm d}$ we write $T_{ij}\in\R$ for the $j$th element in the $i$th row and $I_d\in\R^{d\tm d}$ is the identity matrix. Throughout, $Y_+\subset\R^d$ denotes an order cone inducing the relations $\leq,<$ and $\ll$ (cf.~\eqref{noorders3} and App.~\ref{secA} for related terminology); $Y_+'$ is the dual cone. 

On metric spaces $(X,d)$, $U^\circ$ denotes the interior, $\overline{U}$ the closure of a set $U\subseteq X$, and $B_r(x):=\set{y\in X:\,d(y,x)<r}$ the open ball of radius $r>0$ and $x\in X$. 

With Banach spaces $X,Y$ we write $L(X,Y)$ for the linear space of bounded linear maps $T:X\to Y$, $L(X):=L(X,X)$ and $X':=L(X,\R)$ for the dual space, as well as $N(T):=T^{-1}(\set{0})$ for the kernel of $T$. 
\section{Fredholm integral operators}
\label{sec2}
This section provides sufficient conditions for \emph{Fredholm (integral) operators} 
\begin{equation}
	\tag{$F$}
	\eK u:=\int_\Omega K(\cdot,y)u(y)\d\mu(y)
	\label{nofred}
\end{equation}
to preserve order relations. In essence, we demonstrate that related positivity properties of the kernels $K:\Omega^2\to L(\R^d)$ carry over to integral operators $\eK$. 

For this purpose and referring to later applications it is ambient to work with an abstract measure-theoretical integral in \eqref{nofred}. Thereto, assume $(\Omega,\fA,\mu)$ is a measure space satisfying $\mu(\Omega)<\infty$. The resulting abstract $\mu$-integral of a $\mu$-measurable function $u:\Omega\to\R^d$ is denoted by $\int_\Omega u(y)\d\mu(y)$ and satisfies
\begin{equation}
	\sprod{\int_\Omega u(y)\d\mu(y),y'}
	=
	\int_\Omega\sprod{u(y),y'}\d\mu(y)\fall y'\in\R^d.
	\label{noprod}
\end{equation}
A relevant case in applications \cite{lutscher:19} is the $\kappa$-dimensional Lebesgue measure $\mu=\lambda_\kappa$ on $\Omega\subset\R^\kappa$ yielding the Lebesgue integral $\int_\Omega u(y)\d y:=\int_\Omega u(y)\d\lambda_\kappa(y)$ for functions $u:\Omega\to\R^d$. Moreover, numerical numerical analysis relies on
\begin{rem}[quadrature methods as abstract integrals]\label{remquad}
	Let $\Omega'\subset\R^\kappa$ be countable and $\eta\in\Omega'$ with associated reals $w_\eta\geq 0$. Then
	$
		\mu(\Omega'):=\sum_{\eta\in\Omega'}w_\eta
	$
	is a measure on the family of countable subsets of $\R^\kappa$. Note that $\sum_{\eta\in\Omega'}w_\eta<\infty$ gives $\mu(\Omega')<\infty$ and the $\mu$-integral is $\int_{\Omega'}u(y)\d\mu(y)=\sum_{\eta\in\Omega'}w_\eta u(\eta)$. This setting includes the case, where $\Omega'$ is a finite set $\Omega_n=\set{\eta_j:\,0\leq j\leq q_n}$ of nodes $\eta_j$ with associated weights $w_\eta=w_j$ of a numerical quadrature rule \cite{atkinson:97,hackbusch:95}, where $(q_n)_{n\in\N}$ is a strictly increasing sequence of positive integers. Here, 
	\begin{align*}
		\mu(\Omega_n)&=\sum_{j=0}^{q_n}w_j,&
		\int_{\Omega_n}u(y)\d\mu(y)&=\sum_{j=0}^{q_n}w_ju(\eta_j)
	\end{align*}
	for all functions $u:\Omega_n\to\R^d$. Concrete examples are listed in App.~\ref{secC}. 
\end{rem}
\subsection{Fredholm operators on $\Cd$}
In this subsection, we restrict to measure spaces $(\Omega,\fA,\mu)$ with a compact metric space $\Omega$, a $\sigma$-algebra $\fA$ (containing the Borel sets) and a finite measure $\mu$, i.e.\ $\mu(\Omega)<\infty$. 

The set $\Cd$ of continuous functions $u:\Omega\to\R^d$ is a real Banach space when equipped with the maximum norm $\norm{u}_\infty:=\max_{x\in\Omega}\abs{u(x)}$. Moreover, 
$$
	\Cd_+:=\set{u\in\Cd:\,u(x)\in Y_+\text{ for all }x\in\Omega}
$$
abbreviates the set of continuous functions having values in a cone $Y_+\subset\R^d$. 
\begin{lem}
	The set $\Cd_+$ is a cone. If $Y_+$ is solid, then $\Cd_+$ is solid and total. 
\end{lem}
\begin{pf}
	The closedness and convexity of $Y_+$ extend to $\Cd_+$ and a pointwise consideration yields $\R_+\Cd_+\subseteq\Cd_+$ and $\Cd_+\cap(-\Cd_+)=\set{0}$. Thus, $\Cd_+$ is a cone. If $Y_+\subset\R^d$ is solid, then there exist $y_0\in Y_+$ and $\eps>0$ so that $B_\eps(y_0)\subset Y_+$. Hence, $u_0(x):\equiv y_0$ on $\Omega$ is an interior point of $\Cd_+$ and $(\Cd_+)^\circ\neq\emptyset$ follows. 
	\qed	
\end{pf}

Having identified $\Cd_+$ as (solid) cone, we introduce the relations (cf.~\eqref{noorders3}) 
\begin{align*}
	u\le\bar u\quad&:\Leftrightarrow\quad\bar u-u\in\Cd_+,\\
	u\lle\bar u\quad&:\Leftrightarrow\quad\bar u-u\in\Cd_+\setminus\set{0},\\
	u\lld\bar u\quad&:\Leftrightarrow\quad\bar u-u\in(\Cd_+)^\circ\fall u,\bar u\in\Cd,
\end{align*}
allowing the subsequent characterization:
\begin{lem}\label{lemcone}
	The following holds for $u,\bar u\in\Cd$: 
	\begin{enumerate}
		\item[(a)] $u\le\bar u\Leftrightarrow u(x)\leq\bar u(x)$ for all $x\in\Omega$ $\Leftrightarrow\sprod{u(x),y'}\leq\sprod{\bar u(x),y'}$ for all $x\in\Omega$ and $y'\in Y_+'$. 

		\item[(b)] $u\lle\bar u\Leftrightarrow u(x)\leq\bar u(x)$ for all $x\in\Omega$ and $u(x_0)<\bar u(x_0)$ for some $x_0\in\Omega$.

		\item[(c)] If $Y_+$ is solid, then $u\lld\bar u\Leftrightarrow u(x)\ll\bar u(x)$ for all $x\in\Omega$ $\Leftrightarrow\sprod{u(x),y'}<\sprod{\bar u(x),y'}$ for all $x\in\Omega$, $y'\in Y_+'\setminus\set{0}$.
	\end{enumerate}
\end{lem}
\begin{pf}
	(a) follows immediately by definition and from \lref{lemA1}(a). 
	\\
	(b) results from the definition. 
	\\
	(c) Due to linearity it suffices to establish the claim for the pair $(0,u)$ rather than the functions $(u,\bar u)$. We have to show two directions:\\
	$(\Leftarrow)$ Assume that $0<\sprod{u(x),y'}$ for all $x\in\Omega$, $y'\in Y_+'\setminus\set{0}$ holds. Therefore, \lref{lemA1}(b) ensures that for any $x\in\Omega$ there exists an $\eps_x>0$ guaranteeing the inclusion $B_{\eps_x}(u(x))\subset Y_+$. We prove that
	\begin{equation}\label{ball}
		\exists\eps>0:\forall x\in\Omega:\,B_\eps(u(x))\subset Y_+.
	\end{equation}
	Assuming the contrary there exist sequences $(x_n)_{n\in\N}$ in $\Omega$ and $(z_n)_{n\in\N}$ in $\R^d$ with
	\begin{equation}\label{ball1}
		\abs{z_n-u(x_n))}\leq\tfrac{1}{n}, \quad z_n\notin Y_+\fall n\in\N.
	\end{equation}
	From compactness of $\Omega$ and continuity of $u$ there exists a subsequence $(x_{k_n})_{n\in\N}$ with limit $x^\ast\in\Omega$ such that $\lim_{n\to\infty}z_{k_n}=u(x^\ast)$ is satisfied. Consequently, $z_{k_n}\in B_{\eps_{x^\ast}}(u(x^\ast))\subset Y_+$ holds for sufficiently large $n\in\N$, contradicting \eqref{ball1}. Hence, \eqref{ball} holds and implies that $B_\eps(u)\subset\Cd_+$ and $0\lld u$. 
	\\
	$(\Rightarrow)$ For $0\lld u$ there exists an $\eps>0$ so that $v(x)\in Y_+$ for all $v\in B_\eps(u)$ and $x\in\Omega$. Let $x\in\Omega$ and $z_0\in B_\eps(u(x))\subseteq\R^d$. If $\bar v(\xi):=u(\xi)+z_0-u(x)$, $\xi\in\Omega$, then $\bar v\in B_\eps(u)$ and thus $\bar v(\xi)\in Y_+$ for all $\xi\in\Omega$. For $\xi=x$ we see that $z_0\in Y_+$ and $B_\eps(u(x))\subset Y_+$ for each $x\in\Omega$. Hence, $u(x)\in Y^\circ_+$ for all $x\in\Omega$. From \lref{lemA1}(b) we have $0<\sprod{u(x),y'}$ for all $x\in\Omega$, $y'\in Y_+'\setminus\set{0}$. 
	\qed
\end{pf}

\begin{exmp}\label{exgencone}
	If a set $\set{e_1,\ldots,e_d}\subset\R^d$ is linearly independent, then the linear combinations
	$
		Y_+=\set{\sum_{i=1}^d\alpha_ie_i:\,\alpha_1,\ldots,\alpha_d\geq 0}\subset\R^d
	$
	form a solid, thus total cone with interior $Y_+^\circ=\bigl\{\sum_{i=1}^d\alpha_ie_i:\,\alpha_1,\ldots,\alpha_d>0\bigr\}$. With $\set{e_1',\ldots,e_d'}\subset\R^d$ chosen according to $\sprod{e_i,e_j'}=\delta_{ij}$ for $1\leq i,j\leq d$, then
	$
		Y_+'=\bigl\{\sum_{j=1}^d\beta_je_j':\,\beta_1,\ldots,\beta_d\geq 0\bigr\}\subset\R^d
	$
	is the dual cone to $Y_+$. Given $u,\bar u\in\Cd$ one has:
	\begin{align*}
		u\le\bar u\quad\Leftrightarrow&\quad\sprod{u(x),e_j'}\leq\sprod{\bar u(x),e_j'}
		\text{ for all }x\in\Omega,\,1\leq j\leq d,\\
		u\lle\bar u\quad\Leftrightarrow&\quad\sprod{u(x),e_j'}\leq\sprod{\bar u(x),e_j'}
		\text{ for all }x\in\Omega,\,1\leq j\leq d\\
		&\quad
		\text{and }u(x_0)<\bar u(x_0)\text{ for some }x_0\in\Omega,\\
		u\lld\bar u\quad\Leftrightarrow&\quad\sprod{u(x),e_j'}<\sprod{\bar u(x),e_j'}\text{ for all }x\in\Omega,\,1\leq j\leq d.
	\end{align*}
\end{exmp}
\begin{exmp}[orthants]\label{exorthants}
	Let $Y_+$ be an orthant of $\R^d$ spanned by the linearly independent vectors $e_i=\varsigma_i(\delta_{ij})_{j=1}^d$ with $\varsigma_i=\pm 1$. For $e_j':=\varsigma_j(\delta_{ij})_{i=1}^d$ one obtains characterizations based on the component functions
	\begin{align*}
		u\le\bar u\quad\Leftrightarrow&\quad\varsigma_iu_i(x)\leq\varsigma_i\bar u_i(x)\fall x\in\Omega,\,1\leq i\leq d,\\
		u\lle\bar u\quad\Leftrightarrow&\quad\varsigma_iu_i(x)\leq\varsigma_i\bar u_i(x)\fall x\in\Omega,\,1\leq i\leq d\text{ and there}\\
		&\quad\text{exist } x_0\in\Omega, 1\leq i_0\leq d\text{ with } u_{i_0}(x_0)\neq\bar u_{i_0}(x_0),\\
		u\lld\bar u\quad\Leftrightarrow&\quad\varsigma_iu_i(x)<\varsigma_i\bar u_i(x)\fall x\in\Omega,\,1\leq i\leq d.
	\end{align*}
\end{exmp}

A convenient setting for linear integral operators provides the following
\begin{hypo}
	Assume that $K:\Omega\tm\Omega\to L(\R^d)$ fulfills: 
	\begin{itemize}
		\item[$(L)$] $K(x,\cdot):\Omega\to L(\R^d)$ is $\mu$-measurable for all $x\in\Omega$ with 
		$$
			\sup_{x\in\Omega}\int_\Omega\abs{K(x,y)}\d\mu(y)<\infty
		$$
		and $\lim_{x\to x_0}\int_\Omega\abs{K(x,y)-K(x_0,y)}\d\mu(y)=0$ for all $x_0\in\Omega$. 
	\end{itemize}
\end{hypo}
These assumptions yield that Fredholm operators $\eK$ given in \eqref{nofred} readily fulfill the inclusion $\eK\in L(C(\Omega)^d)$ (see \cite[p.~167, Prop.~3.4]{martin:76}).
In case the limit relation in $(L)$ holds uniformly in $x_0\in\Omega$, then $\eK$ is even compact; however this is not immediately relevant for our further analysis. 
\begin{thm}[positivity of $\eK$ on $\Cd$]\label{thmlinmon}
	Let Hypothesis $(L)$ hold. If a kernel $K(x,y)\in L(\R^d)$ is $Y_+$-positive for all $x\in\Omega$ and $\mu$-a.a.\ $y\in\Omega$, then a Fredholm operator $\eK\in L(\Cd)$ is $\Cd_+$-positive. 
	These assumptions additionally yield:
	\begin{enumerate}
		\item[(a)] If nonempty, open subsets of $\Omega$ have positive measure, there exists a $\bar x\in\Omega$ so that $K(\bar x,\cdot)$ is continuous on $\Omega$ and $K(\bar x,y)$ is $Y_+$-injective for $\mu$-a.a.\ $y\in\Omega$, then $\eK$ strictly $\Cd_+$-positive. 
		
		\item[(b)] If nonempty, open subsets of $\Omega$ have positive measure, $Y_+$ is solid and $K(x,y)$ is strongly $Y_+$-positive for all $x\in\Omega$ and $\mu$-a.a.\ $y\in\Omega$, then $\eK$ is strongly $\Cd_+$-positive.

		\item[(c)] If $\mu(\Omega)>0$, $Y_+$ is solid and $K(x,y)Y_+^\circ\subseteq Y_+^\circ$ for all $x\in\Omega$ and $\mu$-a.a.\ $y\in\Omega$, then $\eK(\Cd_+)^\circ\subseteq(\Cd_+)^\circ$. 
	\end{enumerate}
\end{thm}
\begin{pf}
	Let $x\in\Omega$, $y'\in Y_+'$ and $0\lle u$. By \lref{lemcone}(a) the positivity of $\eK$ is an immediate consequence of
	$$
		\sprod{(\eK u)(x),y'}
		\stackrel{\eqref{noprod}}{=}
		\int_{\Omega}\sprod{K(x,y)u(y),y'}\d\mu(y)\geq 0.
	$$
	Thus, for the remaining proof $\eK$ is positive.
	\\
	(a) Let $0\lle u$. There exists a $y_0\in\Omega$ such that $\max_{x\in\Omega}\abs{u(x)}=\abs{u(y_0)}>0$ and thus $u_{i_0}(y_0)\neq 0$ for some index $i_0\in\{1,\ldots,d\}$. First, if $u_{i_0}(y_0)>0$, then
	$
		\Omega_0
		:=
		\{y\in\Omega: u_{i_0}(y)>\tfrac{1}{2}u_{i_0}(y_0)\}
		=
		u_{i_0}^{-1}((\tfrac{1}{2}u_{i_0}(y_0),\infty))\neq\emptyset
	$
	and the continuity of $u_{i_0}$ yields that $\Omega_0$ is open; thus, $\mu(\Omega_0)>0$ by assumption. Thanks to $Y_+$-injectivity, $K(\bar x, y)u(y)\in Y_+\setminus\{0\}$ holds for $\mu$-a.a.\ $y\in\Omega_0$ and \lref{lemA1}(a) ensures that there exists a functional $y_y'\in Y_+'\setminus\{0\}$ satisfying $\sprod{K(\bar x,y)u(y),y_y'}>0$. Due to the continuity of $K(\bar x,\cdot)$, for $\mu$-a.a.\ $y\in\Omega_0$ there exists $\eps_y>0$ such that
	$
		\sprod{K(\bar x,\eta)u(\eta),y_y'}>0\quad\text{for}\ \eta\in B_{\eps_y}(y).
	$
	Since $\{B_{\eps_y}(y):\,\mu\textrm{-a.a.}\ y\in\Omega_0\}$ is an open cover of the compact closure $\overline{\Omega_1}$, where $\Omega_1$ is an open set and $\overline{\Omega_1}\subset \Omega_0$ the Borel-Lebesgue Theorem yields a finite subcover $\{B_{\eps_i}(y_i):\,1\leq i\leq n\}$ of $\overline{\Omega_1}$. Moreover $\mu(\Omega_1)>0$. If we define $\tilde y':=\sum^n_{i=1}y'_{y_i}\in Y_+'\neq 0$, then $\sprod{K(\bar x,y)u(y),\tilde y'}>0$ for $\mu$-a.a.\ $y\in\Omega_1$ and
	\begin{align*}
		\sprod{(\eK u)(\bar x),\tilde y'}
		&\stackrel{\eqref{nofred}}{=}
		\sprod{\int_\Omega K(\bar x,y)u(y)\d\mu(y),\tilde y'}
		\ge
		\int_{\Omega_1}\sprod{K(\bar x,y)u(y),\tilde y'}\d\mu(y)>0
	\end{align*}
	holds. This implies that $(\eK u)(\bar x)\neq 0$ and therefore $\eK u\neq 0$ holds. Second, in the dual case $u_{i_0}(y_0)<0$ we use the open set $\{y\in\Omega:\, u_{i_0}(y)<\tfrac{1}{2}u_{i_0}(y_0)\}$ instead of $\Omega_0$ and again obtain $\eK u\neq 0$. In conclusion, $\eK$ is strictly positive. 
	\\
	(b) Let $0\lle u$, $y'\in Y_+'\setminus\set{0}$, $x\in\Omega$. Note that $\Omega_0=\{y\in\Omega: u(y)\neq0\}$ is open subset of $\Omega$ and $\mu(\Omega_0)>0$. Since $K(x,y)$ is strongly positive for $\mu$-a.a.\ $y\in\Omega_0$,
	we get $0\ll K(x,y)u(y)$ and hence \lref{lemA1}(b) implies $0<\sprod{K(x,y)u(y),y'}=:\phi(y)$ for $\mu$-a.a.\ $y\in\Omega_0$. Now at least one of the preimages $\Omega_l:=\phi^{-1}((\tfrac{1}{l},\infty))$, $l\in\N$, has positive measure, since $\Omega_0=\bigcup_{l\in\N}\Omega_l$ is of positive measure. In case $\mu(\Omega_l)>0$ we obtain
	$
		0
		<
		\frac{\mu(\Omega_l)}{l}
		\leq
		\int_{\Omega_l}\phi(y)\d\mu(y)
		\leq
		\int_{\Omega_0}\phi(y)\d\mu(y)
	$
	and consequently
	$$
		\sprod{(\eK u)(x),y'}
		\stackrel{\eqref{nofred}}{=}
		\int_{\Omega}\phi(y)\d\mu(y)\ge \int_{\Omega_0}\phi(y)\d\mu(y)\ge \int_{\Omega_l}\phi(y)\d\mu(y)>0.
	$$
	Because $y\in Y_+'\setminus\set{0}$ and $x\in\Omega$ were arbitrary, \lref{lemA1}(b) readily implies $0\ll(\eK u)(x)$ for any $x\in\Omega$ and therefore $\eK$ is strongly $\Cd_+$-positive by \lref{lemcone}(c). 
	\\
	(c) Let $u\in(\Cd_+)^\circ$. Then \lref{lemcone}(c) yields $0\ll K(x,y)u(y)$ for $\mu$-a.a.\ $y\in\Omega$ and $x\in\Omega$ and the claim results as in (b). 
	\qed
\end{pf}
\subsection{Fredholm operators on $L^p(\Omega)$}
Let $p\in[1,\infty)$ and assume $\mu(\Omega)>0$ throughout. The set $L^p(\Omega)^d$ of $p$-in\-te\-grable functions $u:\Omega\to\R^d$ defines a real Banach space with
$$
	\norm{u}_p
	:=
	\intoo{\int_\Omega\abs{u(y)}^p\d\mu(y)}^{1/p}
$$
as norm. Let us furthermore write
$$
	L^p(\Omega)_+^d:=\set{u\in L^p(\Omega)^d:\,u(x)\in Y_+\text{ for $\mu$-a.a.\ }x\in\Omega}
$$
for the set of $p$-integrable functions with values in $Y_+$. 
\begin{lem}
	The set $L^p(\Omega)_+^d$ is a cone. If $Y_+\subseteq\R^d$ is an orthant, then $L^p(\Omega)_+^d$ is total.
\end{lem}
Note that even in case $Y_+=\R_+$ the cone $L^p(\Omega)_+$ is not solid (cf.\ \cite[Ex.~1.11]{amann:76}). Yet, for $Y_+=\R_+$ it is reproducing (i.e.\ $X=Y_+-Y_+$) and hence total. 
\begin{pf}
	Closedness and convexity of $Y_+$ extend to $L^p(\Omega)_+^d$ and a pointwise consideration yields $\R_+L^p(\Omega)_+^d\subseteq L^p(\Omega)_+^d$ and $L^p(\Omega)_+^d\cap(-L^p(\Omega)_+^d)=\set{0}$. So, $L^p(\Omega)_+^d$ is a cone. If $Y_+\subseteq\R^d$ is an orthant, then $L^p(\Omega)_+^d$ is reproducing and consequently total. 
	\qed
\end{pf}
Having identified the set $L^p(\Omega)_+^d$ as a cone, we introduce the relations (cf.~\eqref{noorders3})
\begin{align*}
	u\le\bar u\quad&:\Leftrightarrow\quad\bar u-u\in L^p(\Omega)_+^d,\\
	u\lle\bar u\quad&:\Leftrightarrow\quad\bar u-u\in L^p(\Omega)_+^d\setminus\set{0},
\end{align*}
allowing the subsequent characterization: 
\begin{lem}\label{lemconelp}
	The following holds for $u,\bar u\in L^p(\Omega)^d$: 
	\begin{enumerate}
		\item[(a)] $u\le\bar u\Leftrightarrow u(x)\leq\bar u(x)$ for $\mu$-a.a.\ $x\in\Omega$ $\Leftrightarrow\sprod{u(x),y'}\leq\sprod{\bar u(x),y'}$ for $\mu$-a.a.\ $x\in\Omega$ and all $y'\in Y_+'$. 

		\item[(b)] $u\lle\bar u\Leftrightarrow u(x)\leq\bar u(x)$ $\mu$-a.a.\ in $\Omega$ and $\mu(\set{x\in\Omega:\,u(x)\neq\bar u(x)})>0$. 
	\end{enumerate}
\end{lem}
\begin{pf}
	(a) follows immediately from \lref{lemA1}(a) and by definition. 
	\\
	(b) results from the definition. 
	\qed
\end{pf}

\begin{lem}[strict monotonicity of the integral]\label{lemlp}
	Let $u,\bar u\in L^p(\Omega)_+^d$ w.r.t.\ an orthant $Y_+\subset\R^d$. If $u\lle\bar u$, then $\int_\Omega u(y)\d\mu(y)<\int_\Omega\bar u(y)\d\mu(y)$. 
\end{lem}
\begin{pf}
	It suffices to show that $u\in L^p(\Omega)_+^d\setminus\set{0}$ implies $\int_\Omega u(y)\d\mu(y)\neq 0$. Let $Y_+$ be spanned by the vectors $e_i=\varsigma_i(\delta_{ij})_{j=1}^d$ with $\varsigma_i=\pm 1$. Since $u\neq 0$, there exists an index $j$ and a set $\Omega_0\subseteq\Omega$ of positive measure with $0<\varsigma_ju_j(x)$ for $x\in\Omega_0$. The sets $\Omega_l:=\set{x\in\Omega_0:\,\tfrac{1}{l}\leq\varsigma_ju_j(x)}$ for $l\in\N$ fulfill $\Omega_0=\bigcup_{l\in\N}\Omega_l$ and $\Omega_l\subseteq\Omega_{l+1}$. Thus, the assumption $\int_{\Omega_0}u(y)\d\mu(y)=0$ leads to
	$
		0
		=
		\int_{\Omega_0}\varsigma_ju_j(y)\d\mu(y)
		\geq
		\int_{\Omega_l}\varsigma_ju_j(y)\d\mu(y)
		=
		\tfrac{1}{l}\mu(\Omega_l)
	$
	and hence $\mu(\Omega_l)=0$ for $l\in\N$. Now the continuity of measures implies the contradiction $\mu(\Omega_0)=0$. 
	\qed
\end{pf}

Now we proceed to linear integral operators on the square-integrable functions: 
\begin{hypo}
	Assume that $K:\Omega\tm\Omega\to L(\R^d)$ fulfills: 
	\begin{itemize}
		\item[$(L^2)$] $K:\Omega\tm\Omega\to L(\R^d)$ is $\mu\otimes\mu$-measurable with
		$$
			\int_\Omega\int_\Omega\abs{K(x,y)}^2\d\mu(y)\d\mu(x)<\infty.
		$$
	\end{itemize}
\end{hypo}
Then the Fredholm operator $\eK\in L(L^2(\Omega)^d)$ defined via \eqref{nofred} is well-defined and also compact (cf.\ \cite[p.~47, Thm.~3.2.7]{hackbusch:95}). 
\begin{thm}[positivity of $\eK$ on $L^2(\Omega)^d$]\label{thml2}
	Let Hypothesis $(L^2)$ hold. If a kernel $K(x,y)\in L(\R^d)$ is $Y_+$-positive for $\mu$-a.a.\ $x,y\in\Omega$, then a Fredholm operator $\eK\in L(L^2(\Omega))$ is $L^2(\Omega)_+^d$-positive. If moreover, $Y_+\subset\R^d$ is an orthant and $K(x,y)$ is strictly $L^2(\Omega)_+$-positive for $\mu$-a.a.\ $x,y\in\Omega$, then $\eK$ is strictly $L^2(\Omega)_+^d$-positive. 
\end{thm}
\begin{pf}
	Proving positivity of $\eK$ is formally identical to the argument in the proof of \tref{thmlinmon}. One merely replaces the reference to \lref{lemcone}(a) by \lref{lemconelp}(a). The statement on strict positivity results from \lref{lemlp}. 
	\qed
\end{pf}
\section{Nystr\"om methods}
\label{sec3}
In numerics or simulations of Fredholm operators, the involved integrals can be evaluated only approximately. One achieves this by applying discretization methods from the numerical analysis of integral equations. The most natural and popular discretizations of integral operators are based on Nystr\"om methods (see \cite[pp.~100ff]{atkinson:97} or \cite[pp.~128ff]{hackbusch:95}), where one replaces the integrals by integration (quadrature, cubature) rules. 

In this section, we suppose that $\Omega\subset\R^\kappa$ is compact with Lebesgue measure $\lambda_\kappa(\Omega)>0$. For a continuous function $u:\Omega\to\R^d$, consider the representation
\begin{equation}
	\tag{$Q_n$}
	\int_\Omega u(y)\d y=\sum_{j=0}^{q_n}w_ju(\eta_j)+E_n(u)
	\label{quad}
\end{equation}
with a sequence $(q_n)_{n\in\N}$ in $\N$, \emph{nodes} from a finite set $\Omega_n:=\set{\eta_0,\ldots,\eta_{q_n}}\subseteq\Omega$ and \emph{weights} $w_j\in\R$ such that the error term satisfies $\lim_{n\to\infty}E_n(u)=0$. Such schemes are called \emph{convergent} and we refer to App.~\ref{secC} for concrete examples. 

We say that an integration rule \eqref{quad} fulfills the \emph{net condition}, if 
\begin{equation}
	\forall\eps>0:\exists n_0\in\N:\,\Omega\subseteq\bigcup_{j=0}^{q_n}B_\eps(\eta_j)\fall n\geq n_0(\eps). 
	\label{nonet}
\end{equation}
This assumption is indeed frequently met: 
\begin{exmp}[net condition]
	Condition \eqref{nonet} is satisfied, if the distance between neighboring nodes in $\Omega_n$ can be made arbitrarily small as $n\to\infty$. Therefore, essentially all relevant classes of quadrature formulas \eqref{quad} fulfill the net condition: Providing the nodes $\eta_j^n$ only over the interval $\Omega=[-1,1]$ for simplicity, then Clenshaw-Curtis ($\eta_j^n=\cos(\tfrac{j-1}{n-1}\pi)$, \cite[p.~86]{davis:rabinowitz:84}), Gau{\ss}-Legendre ($\eta_j^n$ are the zeros of the Legendre polynomials $P_n$, \cite[Thm.~5.1]{jordaan:tookos:09}) or Gau{\ss}-Lobatto ($\eta_j^n$ are the zeros of the derivatives $P_n'$) types do work. In each case, one has the limit $\lim_{n\to\infty}\sup_{j=0}^{q_n-1}(\eta_{j+1}^n-\eta_j^n)=0$. Of course, composite quadrature rules (see e.g.\ \cite[pp.~70ff]{davis:rabinowitz:84}) fulfill the net condition throughout. Finally, product cubature rules obtained from the above quadratures (see \cite[pp.~354ff]{davis:rabinowitz:84}) satisfy \eqref{nonet} as well. 
\end{exmp}

We impose the following
\begin{hypo}
	Assume that a kernel $K:\Omega\tm\Omega\to L(\R^d)$ fulfills: 
	\begin{itemize}
		\item[$(NL)$] $K(\cdot,y):\Omega\to L(\R^d)$ is continuous for all $y\in\Omega$.
	\end{itemize}
\end{hypo}
Note that $(NL)$ is sufficient for Hypothesis $(L)$ and $(L^2)$ keeping our previous results applicable. Referring to \rref{remquad}, this leads to the (spatially) \emph{discrete Fredholm operator}
\begin{equation}
	\eK^nu:=\sum_{j=0}^{q_n}w_jK(\cdot,\eta_j)u(\eta_j).
	\label{nofredn}
\end{equation}
There are two natural choices for the domain of $\eK^n$, namely a spatially continuous one $\Cd$ and the spatially discrete function space
$$
	C(\Omega_n)^d=\set{u:\Omega_n\to\R^d};
$$
both are equipped with the $\max$-norm $\norm{\cdot}_\infty$. In each case, $(NL)$ suffices to obtain that $\eK^n$ is well-defined and continuous. 
\begin{rem}[$\eK^n$ on the domain $C(\Omega_n)^d$]\label{remsame}
	Let an integration rule \eqref{quad} have nonnegative weights. In the setting of \rref{remquad}, the above \tref{thmlinmon} applies for the measure $\mu$ from \rref{remquad} and shows that positivity of $\eK^n\in L(C(\Omega_n)^d)$ or $\eK^n\in L(C(\Omega_n)^d,\Cd)$ holds literally with the assumption ``$\mu$-a.a.\ $y\in\Omega$'' replaced by ``all $y\in\Omega_n$''. 
\end{rem}

On the domain $\Cd$ one cannot expect a discrete Fredholm operator $\eK^n$ to be strictly or strongly positive. This is due to the fact that $\Cd_+\setminus\set{0}$ contains functions $u$ vanishing everywhere except from being positive on arbitrarily small domains disjoint from $\Omega_n$. Hence, they are not captured by the Nystr\"om grid $\Omega_n$, that is, $u|_{\Omega_n}=0$ although $u\neq 0$. Consequently, one has $\eK^nu=0$. 

This requires ambient modifications of our above results captured in \rref{remsame}. 
\begin{thm}[positivity of $\eK^n$ on $\Cd$]\label{thmkny}
	Let Hypothesis $(NL)$ hold. 
	\begin{enumerate}
		\item[(a)] If \eqref{quad}, $n\in\N$, have nonnegative weights and $K(x,\eta)$ is $Y_+$-positive for all $x\in\Omega$, $\eta\in\Omega_n$, then $\eK^n\in L(\Cd)$ is $\Cd_+$-positive. 

		\item[(b)] If \eqref{quad}, $n\in\N$, have positive weights, $Y_+$ is solid and $K(x,\eta)Y_+^\circ\subseteq Y_+^\circ$ for all $x\in\Omega$, $\eta\in\Omega_n$, then $\eK^n(\Cd_+)^\circ\subseteq(\Cd_+)^\circ$. 
	\end{enumerate}	
	In case $\Omega=\overline{\Omega^\circ}$, \eqref{quad} have eventually positive weights and the net condition \eqref{nonet} hold, then for each function $u\in\Cd$ with $0\lle u$ there is a $N\in\N$ such that one has for $n\geq N$: 
	\begin{enumerate}
		\item[(c)] If $K(\bar x,\eta)$ is $Y_+$-injective for one $\bar x\in\Omega$ and all $\eta\in\Omega_n$, then $0\lle\eK^n u$. 

		\item[(d)] If $Y_+$ is solid and $K(x,\eta)$ is strongly $Y_+$-positive for all $x\in\Omega$, $\eta\in\Omega_n$, then $0\lld\eK^nu$.
	\end{enumerate}
\end{thm}
A large number of common integration rules \eqref{quad} have positive weights:
\begin{exmp}[positive weights]
	(1) If the set $\Omega$ is a compact interval, then the following classes of quadrature formulas \eqref{quad} have positive weights: Closed Newton-Cotes with $n\in\set{2,\ldots,8,10}$ nodes, open Newton-Cotes with $n\in\set{1,2,4}$ nodes (see \cite[pp.~120--121 and p.~156, Exer.~41]{ralston:rabinowitz:78}), Clenshaw-Curtis \cite[p.~86]{ralston:rabinowitz:78}, Gau{\ss}-Legendre \cite[p.~105]{ralston:rabinowitz:78} and Gau{\ss}-Lobatto. Also composite versions of these quadrature rules clearly have positive weights as well. 
	\\
	(2) If $\Omega$ is a rectangle in $\R^\kappa$, $\kappa>1$, then the product cubature rules obtained from the above quadrature methods feature positive weights. The same holds for domains $\Omega\subset\R^\kappa$ being $C^1$-diffeomorphic to a rectangle $Q\subset\R^\kappa$ by means a $C^1$-diffeomorphism $T:Q\to\Omega$ due to the change of variables formula
	$$
		\int_\Omega u(y)\d y
		=
		\int_Qu(T(x))\abs{\det DT(x)}\d x
	$$
	and applying a cubature rule over rectangles to the right-hand side integral. 
\end{exmp}
\begin{pf}
	(a) By \tref{thmlinmon} with the measure $\mu$ from \rref{remquad} we know $\eK^n\in L(C(\Omega_n)^d)$ is positive. Referring to the definition \eqref{nofredn} this extends to $\eK^n\in L(\Cd)$.
	\\
	(b) results from a direct computation using \lref{lemcone}(c). 
	\\
	For the remaining proof we suppose $0\lle u$. This implies $0\leq u(x)$ for $x\in\Omega$ and there exists $x_0\in\Omega$ with $0\neq u(x_0)$. Let us write $	B:=\{x\in\Omega:\,u(x)\neq 0\}$. By the continuity of $u$ and the assumption on $\Omega$ there exist $\eps_0>0$, $x_1\in\Omega$ so that $B_{2\eps_0}(x_1)\subset B$. Since \eqref{quad} fulfills the net condition \eqref{nonet}, there exists an $n_0(\eps_0)\in\N$ such that for $n\geq n_0(\eps_0)$ there exists $j_{n}\in\{0,\ldots,q_n\}$ with nodes satisfying $\eta_{j_n}\in B$. In the following, assume the rule \eqref{quad} has positive weights for $n\geq n_1$ and let $n\geq N:=\max\set{n_0(\eps_0),n_1}$: 
	\\
	(c) From $Y_+$-injectivity of $K(\bar x,\eta_{j_n})$ one obtains $0<K(\bar x,\eta_{j_n})u(\eta_{j_n})$. Consequently, $0\neq\eK^n(u)(\bar x)$ and hence $0\lle\eK^nu$. 
	\\
	(d) Let $y'\in Y_+'\setminus\{0\}$ and $x\in\Omega$. Using the strong positivity of $K(x,\eta_{j_n})$ we obtain $0\ll K(x,\eta_{n_j})u(\eta_{n_j})$ and \lref{lemA1}(b) implies $0<\sprod{K(x,\eta_{j_n})u(\eta_{j_n}),y'}$. Thus, $0<\sprod{(\eK^nu)(x),y'}$ and using \lref{lemcone}(c) this means $0\lld\eK^n(u)$.
	\qed
\end{pf}
\begin{rem}[convergence of Nystr\"om methods]\label{remnyst}
	While our focus is the persistence of monotonicity under Nystr\"om discretizations, beyond nonnegative weights, reasonable applications also require convergent methods \eqref{quad}, i.e.\
	$$
		\lim_{n\to\infty}\sum_{j=0}^{q_n}w_ju(\eta_j)
		=
		\int_\Omega u(y)\d y\fall u\in\Cd.
	$$
	Following \cite[p.~39, Thm.~3.43]{hackbusch:14} this implies that the quadrature/cubature weights in \eqref{quad} satisfy $\sup_{n\in\N}\sum_{j=0}^{q_n}w_j<\infty$. 
\end{rem}
\section{Projection methods}
\label{sec4}
Let $X(\Omega)$ denote a normed space of functions $u:\Omega\to\R$ being for instance continuous $X(\Omega)=C(\Omega)$ or square-integrable $X(\Omega)=L^2(\Omega)$. Projections methods approximate elements of an infinite-di\-men\-sio\-nal function space $X(\Omega)$ by elements from suitable finite-dimensional subspaces $X_n$, $n\in\N$. 

For this purpose, we choose linearly independent functions $\phi_1,\ldots,\phi_{d_n}\in X(\Omega)$ and linearly independent $\psi_1',\ldots,\psi_{d_n}'\in X(\Omega)'$. Then $X_n:=\spann\set{\phi_1,\ldots,\phi_{d_n}}$ is denoted as \emph{ansatz space} and has dimension $d_n$. If $\phi_1',\ldots,\phi_{d_n}'\in X(\Omega)'$ are functionals satisfying $\phi_i'(\phi_j)=\delta_{ij}$ for all $1\leq i,j\leq d_n$, then 
\begin{align}
	\pi_n:X(\Omega)&\to X_n,&
	\pi_nu&:=\sum_{j=1}^{d_n}\phi_j'(u)\phi_j
	\label{noproj}
\end{align}
defines a bounded projection onto $X_n$. The above functionals are related by
\begin{equation}
	P_n
	\begin{pmatrix}
		\phi_1'\\
		\vdots\\
		\phi_{d_n}'
	\end{pmatrix}
	=
	\begin{pmatrix}
		\psi_1'\\
		\vdots\\
		\psi_{d_n}'
	\end{pmatrix}
	\label{nolgs}
\end{equation}
with an (invertible) matrix $P_n\in\R^{d_n\tm d_n}$. This approach extends to vector-valued functions $u:\Omega\to\R^d$ and operators acting on a Cartesian product $X(\Omega)^d$ in terms of the projection 
\begin{align}
	\Pi_n&\in L(X(\Omega)^d,X_n^d),&
	\Pi_nu
	&:=
	\begin{pmatrix}
		\pi_n u_1\\
		\vdots\\
		\pi_n u_d
	\end{pmatrix}
	=
	\sum^{d_n}_{j=1}
	\phi_j(\cdot)
	\begin{pmatrix}
		\phi_j'(u_1)\\
		\vdots\\
		\phi_j'(u_d)
	\end{pmatrix}.
	\label{no41s}
\end{align}
\begin{rem}[discrete projection methods]
	Projection methods applied to integral operators $\eK$ in \eqref{nofred} merely yield semi-discretizations, that is, although the operators $\Pi_n\eK$ are finite-dimensional, they still contain integrals to be evaluated. Hence, in order to arrive at full discretizations, it remains to apply a Nystr\"om method to $\Pi_n\eK$. In this procedure, positivity properties are preserved, provided the quadrature rules \eqref{quad} satisfy the criteria derived in Sect.~\ref{sec3}. 
\end{rem}
\begin{rem}[convergence of projection methods]
	In the spirit of the above \rref{remnyst} it is reasonable to restrict to convergent projection methods \eqref{no41s}, that is
	$$
		\lim_{n\to\infty}\Pi_nu=u\fall u\in X(\Omega)^d.
	$$
	Thanks to \cite[p.~175, Exer.~8.23]{hackbusch:14} this in turn yields $\sup_{n\in\N}\norm{\Pi_n}<\infty$. 
\end{rem}
\subsection{Collocation methods}
Let $X(\Omega)=C(\Omega)$ be the space of continuous functions over a compact set $\Omega\subset\R^\kappa$ (see \cite[pp.~50ff]{atkinson:97} or \cite[pp.~81ff]{hackbusch:95}) and equip $\R^d$ with the norm
\begin{equation} 
	\abs{x}:=\max_{i=1}^d\abs{x_i}.
	\label{maxnorm}
\end{equation}
For pairwise different \emph{collocation points} $x_1,\ldots,x_{d_n}\in\Omega$ we require the interpolation conditions $(\pi_nu)(x_i)=u(x_i)$ for $1\leq i\leq d_n$ yielding $\psi_i'(u)=u(x_i)$ and resulting in the \emph{collocation matrix} 
$$
	P_n=(\phi_i(x_j))_{i,j=1}^{d_n}. 
$$
\begin{thm}[positivity of $\Pi_n$ on $\Cd$]\label{thmcol}
	If all the functions
	\begin{align*}
		\sigma_i:\Omega&\to\R,&
		\sigma_i(x)&:=\sum_{j=1}^{d_n}(P_n^{-1})_{ij}\phi_j(x)\fall 1\leq i\leq d_n
	\end{align*}
	have nonnegative values, then the following hold:
	\begin{enumerate}
		\item $\Pi_n$ is $\Cd_+$-positive. 

		\item If additionally $Y_+$ is solid and
		\begin{equation}
			\forall x\in\Omega:\exists i_0\in\{1,\ldots,d_n\}:\,\sigma_{i_0}(x)>0
			\label{monstar}
		\end{equation}
		holds, then $\Pi_n(\Cd_+)^\circ\subset (\Cd_+)^\circ$.
	\end{enumerate}
\end{thm}
\begin{pf}
	(a) Let $u\in\Cd$, $0\lle u$ and $y'\in Y_+'$. Due to \lref{lemcone}(a) this means that $0\leq\sprod{u(x),y'}$ for all $x\in\Omega$. If we briefly write $\phi_j'(u)=(\phi_j'(u_k))_{k=1}^d$, then it results 
	\begin{align}
		&
		\sprod{(\Pi_nu)(x),y'}
		\stackrel{\eqref{no41s}}{=}
		\sprod{\sum_{j=1}^{d_n}\phi_j(x)\phi_j'(u),y'}
		=
		\sum_{j=1}^{d_n}\phi_j(x)\sprod{\phi_j'(u),y'}
		\notag\\
		&\stackrel{\eqref{nolgs}}{=}
		\sum_{j=1}^{d_n}\phi_j(x)\sprod{\sum_{i=1}^{d_n}(P_n^{-1})_{ij}u(x_i),y'}
		=
		\sum_{i=1}^{d_n}\sprod{u(x_i),y'}\sum_{j=1}^{d_n}\phi_j(x)(P_n^{-1})_{ij}
		\notag\\
		&=
		\sum_{i=1}^{d_n}\underbrace{\sprod{u(x_i),y'}}_{\geq 0}\sigma_i(x)\geq 0
		\fall x\in\Omega
		\label{nocol3}
	\end{align}
	by assumption. Therefore, \lref{lemcone}(a) yields the claim. 
	\\
	(b)	Let $u\in\Cd$, $0\lld u$. Due to \lref{lemcone}(c) this means $0<\sprod{u(x),y'}$ for all $x\in\Omega$, $y'\in Y_+'\setminus\{0\}$. The relations \eqref{nocol3} and our assumptions ensure that 
	$$
	\sprod{(\Pi_nu)(x),y'}=\sum_{i=1}^{d_n}\sprod{u(x_i),y'}\sigma_i(x)\geq \sprod{u(x_{i_0}),y'}\sigma_{i_0}(x)> 0
	$$ 
	for all $x\in\Omega$, $y'\in Y_+'\setminus\{0\}$ and thus \lref{lemcone}(c) yields $\Pi_nu\in(\Cd_+)^\circ$. 
	\qed
\end{pf}

On the one hand, the relevance of \tref{thmcol} in obtaining structure-preserving collocation methods manifests as follows:
\begin{cor}
	Let Hypothesis $(L)$ hold and $\sigma_i:\Omega\to\R$, $1\leq i\leq d$, have nonnegative values. 
	\begin{enumerate}
		\item[(a)] If $\eK$ is $\Cd_+$-positive, then $\Pi_n\eK\in L(\Cd,X_n^d)$ and $\eK\Pi_n\in L(\Cd)$ are $\Cd_+$-positive.

		\item[(b)] If $\eK$ is strongly $\Cd_+$-positive and \eqref{monstar} holds, then $\Pi_n\eK\in L(\Cd,X_n^d)$ is strongly $\Cd_+$-positive.
	\end{enumerate}
\end{cor}
\begin{pf}
	Statement (a) results since positivity is preserved under composition, while (b) is a consequence of \cref{corA2}.
	\qed
\end{pf}

On the other hand, the applicability of \tref{thmcol} is hindered by the following fact: Many bases $\set{\phi_1,\ldots,\phi_{d_n}}$ (e.g.\ $B$-splines, Bernstein polynomials, etc.) consist of functions having nonnegative values yielding a nonnegative collocation matrix $P_n$. Thus, the inverse $P_n^{-1}$ has nonnegative entries, if and only if $P_n$ is a \emph{monomial} matrix, i.e.\ every column/row contains exactly one positive element (cf.~\cite[p.~2, Thm.~1.1]{kaczorek:02}). 
\begin{rem}[basis transformation]
	The assumption of \tref{thmcol} is invariant under basis transformations. Indeed, if $T\in\R^{d_n\tm d_n}$ is invertible, then also the functions
	$$
		\bar\phi_i:=\sum_{l=1}^{d_n}T_{il}\phi_l:\Omega\to\R\fall 1\leq i\leq d_n
	$$
	define a basis $\set{\bar\phi_1,\ldots,\bar\phi_{d_n}}$ of $X_n$. Due to $\bar\phi_i(x_j)=\sum_{l=1}^{d_n}T_{il}\phi_l(x_j)$ it has the collocation matrix $\bar P_n=TP_n$. This implies
	\begin{align*}
		\bar\sigma
		&:=
		\begin{pmatrix}
			\bar\sigma_1\\
			\vdots\\
			\bar\sigma_{d_n}
		\end{pmatrix}
		=
		\sum_{j=1}^{d_n}
		\begin{pmatrix}
			(\bar P_n^{-1})_{1j}\bar \phi_j\\
			\vdots\\
			(\bar P_n^{-1})_{d_nj}\bar \phi_j
		\end{pmatrix}
		=
		\bar P_n^{-1}
		\begin{pmatrix}
			\bar \phi_1\\
			\vdots\\
			\bar \phi_{d_n}
		\end{pmatrix}\\
		&=
		\bar P_n^{-1}
		\sum_{l=1}^{d_n}
		\begin{pmatrix}
			T_{1l}\phi_l\\
			\vdots\\
			T_{d_nl}\phi_l
		\end{pmatrix}
		=
		\bar P_n^{-1}T
		\begin{pmatrix}
			\phi_1\\
			\vdots\\
			\phi_{d_n}
		\end{pmatrix}
		=
		P_n^{-1}
		\begin{pmatrix}
			\phi_1\\
			\vdots\\
			\phi_{d_n}
		\end{pmatrix}
	\end{align*}
	and consequently $\bar\sigma_i=\sigma_i$ for all $1\leq i\leq d_n$. 
\end{rem}
\subsubsection{Lagrange bases}
Suppose that $\set{\phi_1,\ldots,\phi_{d_n}}$ is a \emph{Lagrange basis} of the ansatz space $X_n$, i.e.\ it satisfies the interpolation condition $\phi_i(x_j)=\delta_{ij}$ for all $1\leq i,j\leq d_n$. Hence, the collocation matrix $P_n=I_{d_n}$ is the identity and the projections read as
$$
	\Pi_nu=\sum_{j=0}^{d_n}\phi_j(\cdot)u(x_j)\fall u\in\Cd.
$$
Concrete examples are as follows, in which $\Omega=[a,b]$ is equipped with a grid \begin{equation}
	a=:x_0<x_1<\ldots<x_{n-1}<x_n:=b
	\label{grid}
\end{equation}
and we abbreviate $h_j:=x_{j+1}-x_j$, $0\leq j<n$. 
\begin{exmp}[piecewise linear collocation]\label{expw}
	Let us consider the ansatz space $X_n\subset C[a,b]$ of piecewise affine functions with break points $\set{x_0,\ldots,x_n}$. It is of dimension $d_n=n+1$ and the \emph{hat functions} $\phi_0,\ldots,\phi_n:[a,b]\to\R_+$, 
	\begin{align}
		\phi_0(x)&:=
		\begin{cases}
			\tfrac{x_1-x}{h_0},&a\leq x\leq x_1,\\
			0,&\text{else},
		\end{cases}\quad
		\phi_n(x):=
		\begin{cases}
			\tfrac{x-x_{n-1}}{h_{n-1}},&x_{n-1}\leq x\leq b,\\
			0,&\text{else},
		\end{cases}
		\notag\\
		\phi_j(x)&:=
		\begin{cases}
			\tfrac{x-x_{j-1}}{h_{j-1}},&x_{j-1}\leq x\leq x_j,\\
			\tfrac{x_{j+1}-x}{h_j},&x_j<x\leq x_{j+1},\\
			0,&\text{else}
		\end{cases}
		\fall 0<j<n
		\label{hatdef}
	\end{align}
	are a suitable basis yielding a convergent method. This extends to rectangles $\Omega:=[a_1,b_1]\tm\ldots\tm[a_\kappa,b_\kappa]$, where each interval $[a_j,b_j]$ may have $n$ subdivisions given by the break points
	$
		a_j=:x_0^j<x_1^j<\ldots<x_{n-1}^j<x_n^j:=b_j
	$
	for $1\leq j\leq\kappa$. 
	If $\phi_i^j:[a_j,b_j]\to\R_+$, $0\leq i\leq n$, are the corresponding hat functions associate to an interval $[a_j,b_j]$, $1\leq j\leq\kappa$, then we define their multivariate version
	$$
		\phi_\iota(\xi):=\prod_{j=1}^\kappa\phi_{\iota_j}^j(\xi_j)
		\text{ for }
		\xi=(\xi_1,\ldots,\xi_\kappa)\in\Omega,\,\iota=(\iota_1,\ldots,\iota_\kappa)\in\set{0,\ldots,n}^\kappa
	$$
	and choose $\set{\phi_\iota:\Omega\to\R_+\mid\iota\in\set{0,\ldots,n}^\kappa}$ as basis of $X_n\subset C(\Omega)$. It has dimension $d_n=(n+1)^\kappa$. It is not hard to see that $\phi_\iota(x)\in[0,1]$ holds on $\Omega$. This yields a Lagrange basis consisting of nonnegative functions. In each case \tref{thmcol} implies that both $\Pi_n$ is $\Cd_+$-positive and $\Pi_n(\Cd_+)^\circ\subseteq(\Cd_+)^\circ$.
\end{exmp}

\begin{exmp}[polynomial interpolation]\label{exlag}
	Let us consider the ansatz space $X_n:=\spann\set{1,x,\ldots,x^n}$ consisting of all polynomials with degree $\leq n$. It is of dimension $d_n=n+1$ and the \emph{Lagrange functions}
	\begin{align*}
		\phi_j:[a,b]&\to\R,&
		\phi_j(x)&:=\prod_{\stackrel{k=0}{k\neq j}}^n\frac{x-x_k}{x_j-x_k}\fall 0\leq j\leq n
	\end{align*}
	yield a Lagrange basis $\set{\phi_0,\ldots,\phi_n}$ for $X_n$ of functions with positive and negative values. Interpolating the constant function $u(x)\equiv 1$ yields that the basis functions form a partition of unity, i.e.\ $\sum_{j=0}^{d_n}\phi_j(x)\equiv 1$ on $[a,b]$. Thus \tref{thmcol} implies both $\Pi_n$ is $C[a,b]_+^d$-positive and $\Pi_n(C[a,b]_+^d)^\circ\subset(C[a,b]_+^d)^\circ$.
\end{exmp}
Although polynomial interpolation preserves positivity, it is burdened by its numerical instability and thus its lack of convergence \cite[p.~84]{hackbusch:95} or \cite[pp.~51--52, Sect.~4.5]{hackbusch:14}; see also \eref{exlaplace}. We next illustrate that various typically used collocation methods do not yield order preserving schemes. 
\subsubsection{Spline interpolation}
\label{sec412}
\begin{figure}
	\includegraphics[scale=0.35]{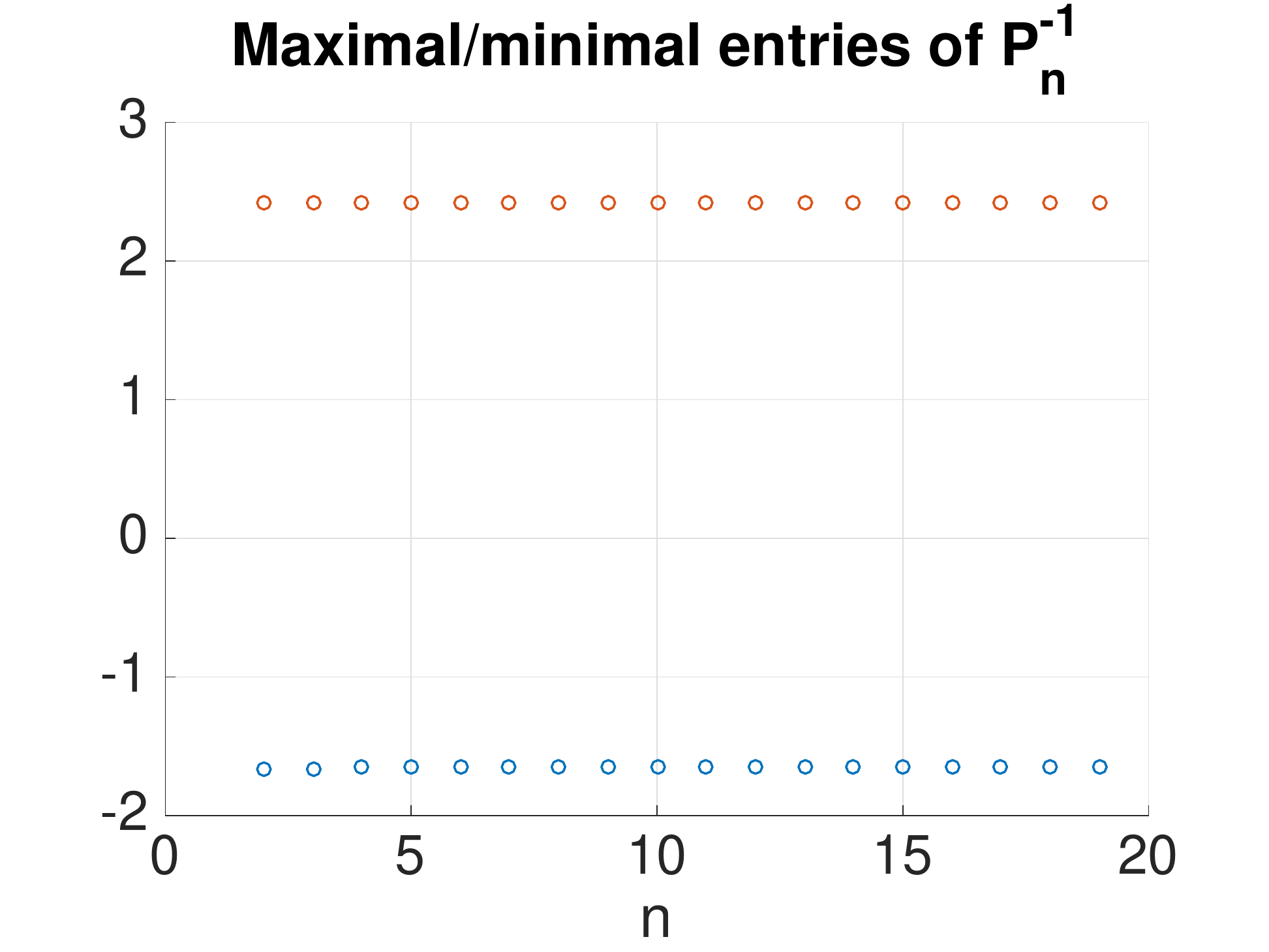}
	\includegraphics[scale=0.35]{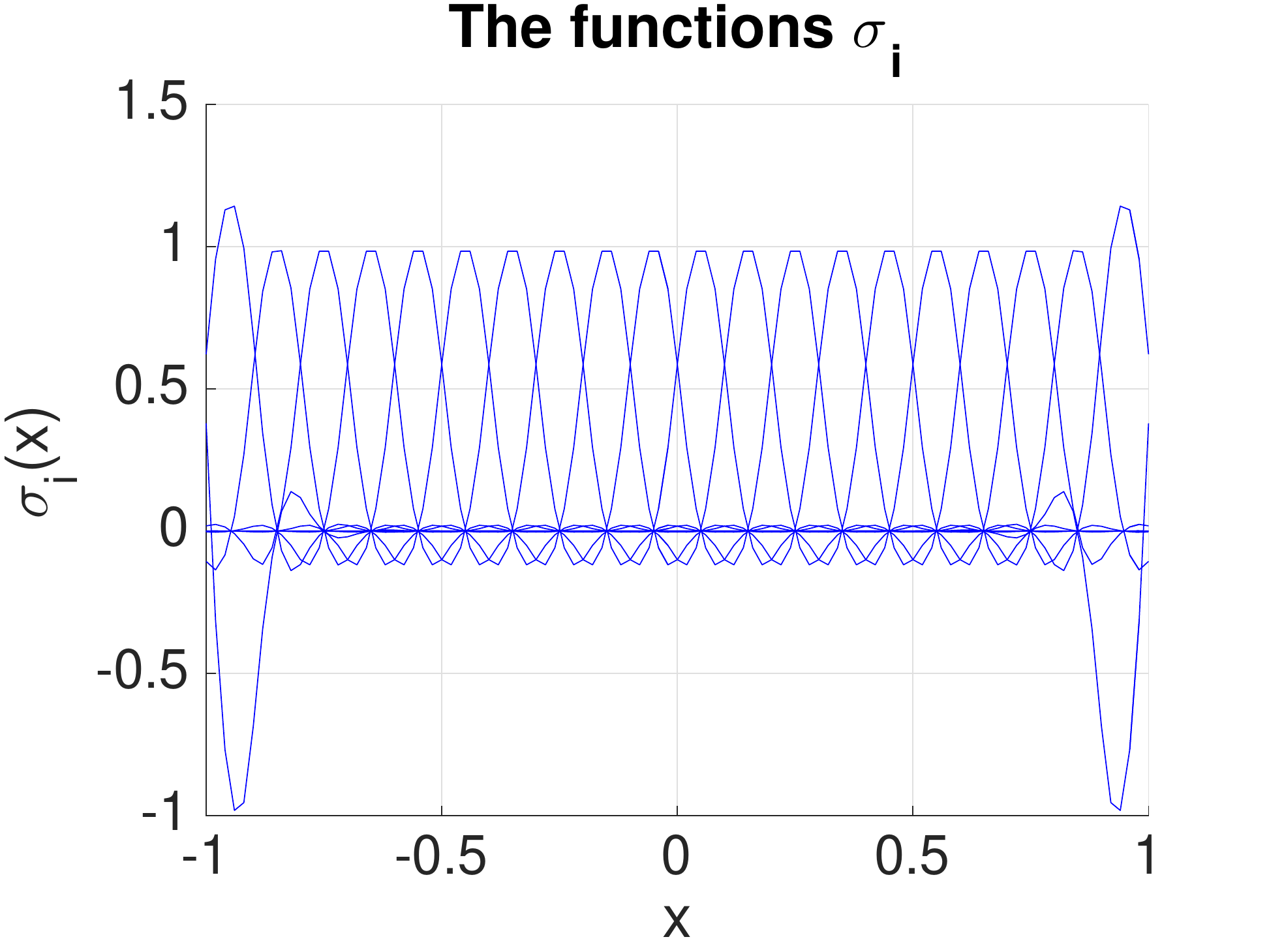}
	\caption{Extremal entries of the inverse collocation matrices $P_n^{-1}$ for the quadratic splines from Sect.~\ref{sec412} depending on $n$ (left) and the associate functions $\sigma_i:[-1,1]\to\R$, $0\leq i\leq n+1$, for $n=20$ (right)}
	\label{figquadsplines}
\end{figure}

\paragraph*{Quadratic splines}
Define $\xi_i:=a+ih$ for $0\leq i\leq n$ with $h:=\tfrac{b-a}{n}$. Let $X_n\subseteq C^1[a,b]$ denote the $d_n=n+2$-dimensional space of quadratic splines equipped with the basis $\set{\phi_{-2},\ldots,\phi_{n-1}}$ of $B$-splines (cf.~\cite[pp.~242ff]{hammerlin:hoffmann:91}) 
\begin{equation}
	\phi_j(x)
	:=
	\frac{1}{2h^2}
	\begin{cases}
		(x-\xi_j)^2,&\xi_j\leq x<\xi_{j+1},\\
		h^2+2h(x-\xi_{j+1})-2(x-\xi_{j+1})^2,&\xi_{j+1}\leq x<\xi_{j+2},\\
		(\xi_{j+3}-x)^2,&\xi_{j+2}\leq x<\xi_{j+3};
	\end{cases}
	\label{phidef}
\end{equation}
they are nonnegative and satisfy $\supp\phi_j=[\xi_j,\xi_{j+3}]$. For the collocation points
\begin{align*}
	x_0&:=a,&
	x_i&:=\tfrac{\xi_i+\xi_{i-1}}{2}\fall 1\leq i\leq n,&
	x_{n+1}&:=b
\end{align*}
establishing the grid \eqref{grid} the collocation matrix becomes tridiagonal
$$
	P_n=
	\frac{1}{2}
	\begin{pmatrix}
		1 & 1 & & &\\
		\tfrac{1}{4} & \tfrac{3}{2} & \tfrac{1}{4} & &\\
		& \ddots & \ddots & \ddots & \\
		&& \tfrac{1}{4} & \tfrac{3}{2} & \tfrac{1}{4} \\
		&&& 1 & 1
	\end{pmatrix}\in\R^{(n+2)\tm(n+2)}
$$
(see \cite[pp.~270]{hammerlin:hoffmann:91}). The extremal entries of $P_n^{-1}$ depending on $n\geq 3$ are illustrated in \fref{figquadsplines}. They have positive and negative signs, and hence \tref{thmcol} does not apply. In order to demonstrate that positivity of the projection operator is indeed violated consider $\xi:=a+h$. One obtains
$$
	(\pi_nu)(\xi)
	\stackrel{\eqref{noproj}}{=}
	\sum_{j=-2}^{n-1}\phi_j'(u)\phi_j(\xi)
	=
	\phi_{-1}'(u)\phi_{-1}(\xi)+\phi_0'(u)\phi_0(\xi)
	\stackrel{\eqref{phidef}}{=}
	\tfrac{\phi_{-1}'(u)}{2}+\tfrac{\phi_0'(u)}{2}
$$
and in particular for $u(x_j)=0$ with $j>0$ one has
$$
	(\pi_nu)(\xi)
	=
	\tfrac{\phi_{-1}'(u)+\phi_0'(u)}{2}
	\stackrel{\eqref{nolgs}}{=}
	\tfrac{(P_n^{-1})_{21}+(P_n^{-1})_{31}}{2}u(a), 
$$
where $u:[a,b]\to\R_+$ is any continuous function with $u(a)>0$ and $u(x)=0$ for $x\geq x_1$. Now \lref{leminv} yields that $(P_n^{-1})_{21}+(P_n^{-1})_{31}<0$ and therefore $\pi_nu\not\in C[a,b]_+$ holds. 
\paragraph*{Cubic splines}
In contrast, let us next construct positivity preserving projection methods as follows: 
\begin{cor}[cubic splines]\label{corspline}
	The piecewise defined cubic spline 
	\begin{align}
		(\Pi_nu)(x)
		:=&
		\intoo{1-3(\tfrac{x-x_i}{h_i})^2+2(\tfrac{x-x_i}{h_i})^3}u(x_i)
		\label{corspline1}\\
		&+
		\intoo{3(\tfrac{x-x_i}{h_i})^2-2(\tfrac{x-x_i}{h_i})^3}u(x_{i+1})
		\fall x\in[x_i,x_{i+1}],\,0\leq i<n
		\notag
	\end{align}
	yields a $C[a,b]_+^d$-positive projection $\Pi_n:C[a,b]^d\to C^1[a,b]^d$ satisfying the estimate $\norm{\Pi_n}_{L(C[a,b]^d)}\leq 2$. 
\end{cor}
Although \eqref{corspline1} gives rise to a positivity preserving spline, its derivative vanishes in all the collocation points $x_j$ and thus $\Pi_nu$ has rather unpleasant approximation properties. This is also exemplified by a low convergence rate in numerical simulations (see \eref{exlaplace}) based on \eqref{corspline1}. 
\begin{pf}
	(I) Given nonnegative real numbers $s_0,\ldots,s_n$, the cubic $C^1$-spline $s:[0,1]\to\R$, piecewise defined as
	\begin{align*}
		s(x)
		=
		s_i+(x-x_i)m_i
		&+(\tfrac{x-x_i}{h_i})^2\intoo{3(s_{i+1}-s_i)-h_i(2m_i+m_{i+1})}\\
		&+(\tfrac{x-x_i}{h_i})^3\intoo{h_i(m_i+m_{i+1})-2(s_{i+1}-s_i)}
	\end{align*}
	for $x\in[x_i,x_{i+1}]$, $0\leq i<n$ is investigated in \cite{schmidt:hess:88}. It satisfies the interpolation conditions $s(x_i)=s_i$, as well as $s'(x_{i})=m_{i}$ for $0\leq i\leq n$, where $m_i\in\R$ are free parameters. According to \cite[Thm.~4]{schmidt:hess:88}, the spline $s$ is nonnegative on the interval $[0,1]$ if and only if $(m_i,m_{i+1})\in W_i(u)$ holds for all $0\leq i<n$, where
	\begin{align*}
		W_i:=&\set{(x,y)\in\R^2:\,-3s_i\leq h_{i+1}x,\,h_{i+1}y\leq 3s_{i+1}}\\
		&\cup
		\{(x,y)\in\R^2:\,0\leq 36s_is_{i+1}(x^2+xy+y^2-3\tau_{i+1}(x+y)+3\tau_{i+1}^2)\\
		&+3(s_{i+1}x-s_iy)(2h_{i+1}xy-3s_{i+1}x+3s_iy)\\
		&+4h_{i+1}(s_{i+1}x^3-s_iy^3)-h_{i+1}^2x^2y^2\}
	\end{align*}
	with $\tau_{i+1}:=\tfrac{s_{i+1}-s_i}{h_i}$. In particular, choosing $m_i:=0$ for all $0\leq i\leq n$ yields the inclusion $(0,0)\in W_i$ and a nonnegative spline
	$$
		s(x)
		=
		(1-3\theta^2+2\theta^3)s_i+(3\theta^2-2\theta^3)s_{i+1}
		\fall x\in[x_j,x_{j+1}],
	$$
	where $\theta:=\tfrac{x-x_i}{h_i}\in[0,1]$. 
	\\
	(II) Let $u\in\Cd_+$, which by \lref{lemcone}(a) means that $0\leq\iprod{u(x),y'}$ for all $x\in[a,b]$, $y'\in Y_+'$. Then the nonnegativity of the spline from step (I) yields
			$$
		\iprod{(\Pi_nu)(x),y'}
		\stackrel{\eqref{corspline1}}{=}
		(1-3\theta^2+2\theta^3)\iprod{u(x_i),y'}+(3\theta^2-2\theta^3)\iprod{u(x_{i+1}),y'}
		\geq 
		0
	$$
	for all $x\in[x_i,x_{i+1}]$, $0\leq i<n$. \lref{lemcone}(a) guarantees $\Pi_nu\in\Cd_+$, that is the positivity claim for $\Pi_n$. Finally, the estimate
	$$
		(\pi_nu)(x)
		=
		(\underbrace{1-3\theta^2+2\theta^3}_{\in[0,1]})u(x_i)+(\underbrace{3\theta^2-2\theta^3}_{\in[0,1]})u(x_{i+1})
		\fall x\in[x_i,x_{i+1}]
	$$
	implies $\abs{(\pi_n u)(x)}\leq 2\norm{u}_\infty$ for all $x\in[a,b]$ and $u\in C[a,b]$. Keeping an eye on \eqref{maxnorm} this leads to the claimed norm estimate for the projections $\Pi_n$. 
	\qed
\end{pf}

\begin{rem}[positivity preserving splines]
	The above \cref{corspline} is based on the idea to apply positivity preserving interpolation methods. Besides the cubic $C^1$-splines from \cite{schmidt:hess:88}, this can alternatively be done using the rational splines of order $\geq 3$ constructed in \cite{sakai:schmidt:89,zhu:18}, where the latter reference provides even error estimates. Via bivariate splines this approach also extends to higher-dimensional domains $\Omega$. \end{rem}
\subsubsection{$\text{sinc}$-Collocation}
\label{exsinc}
\begin{figure}
	\includegraphics[scale=0.35]{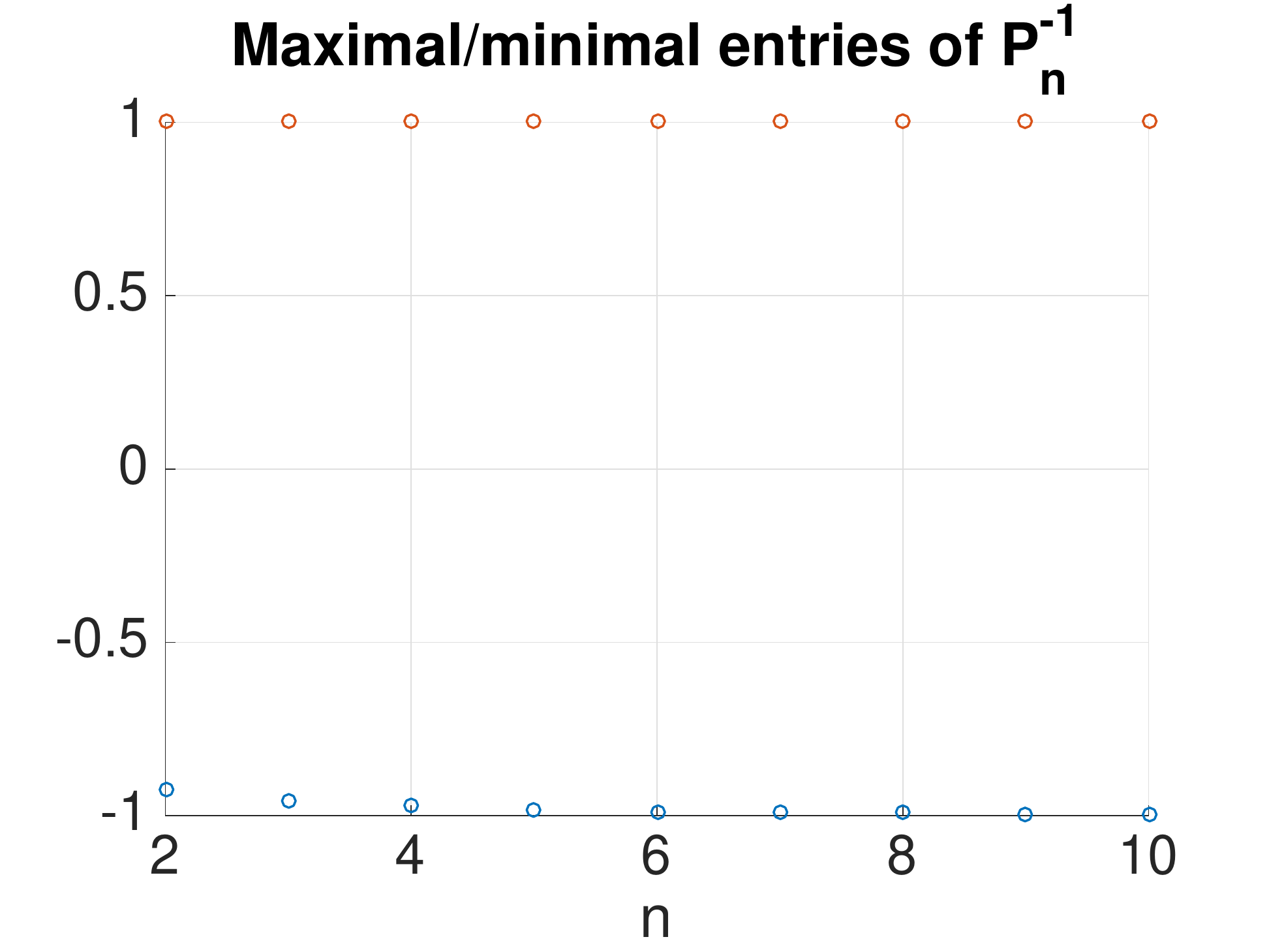}
	\includegraphics[scale=0.35]{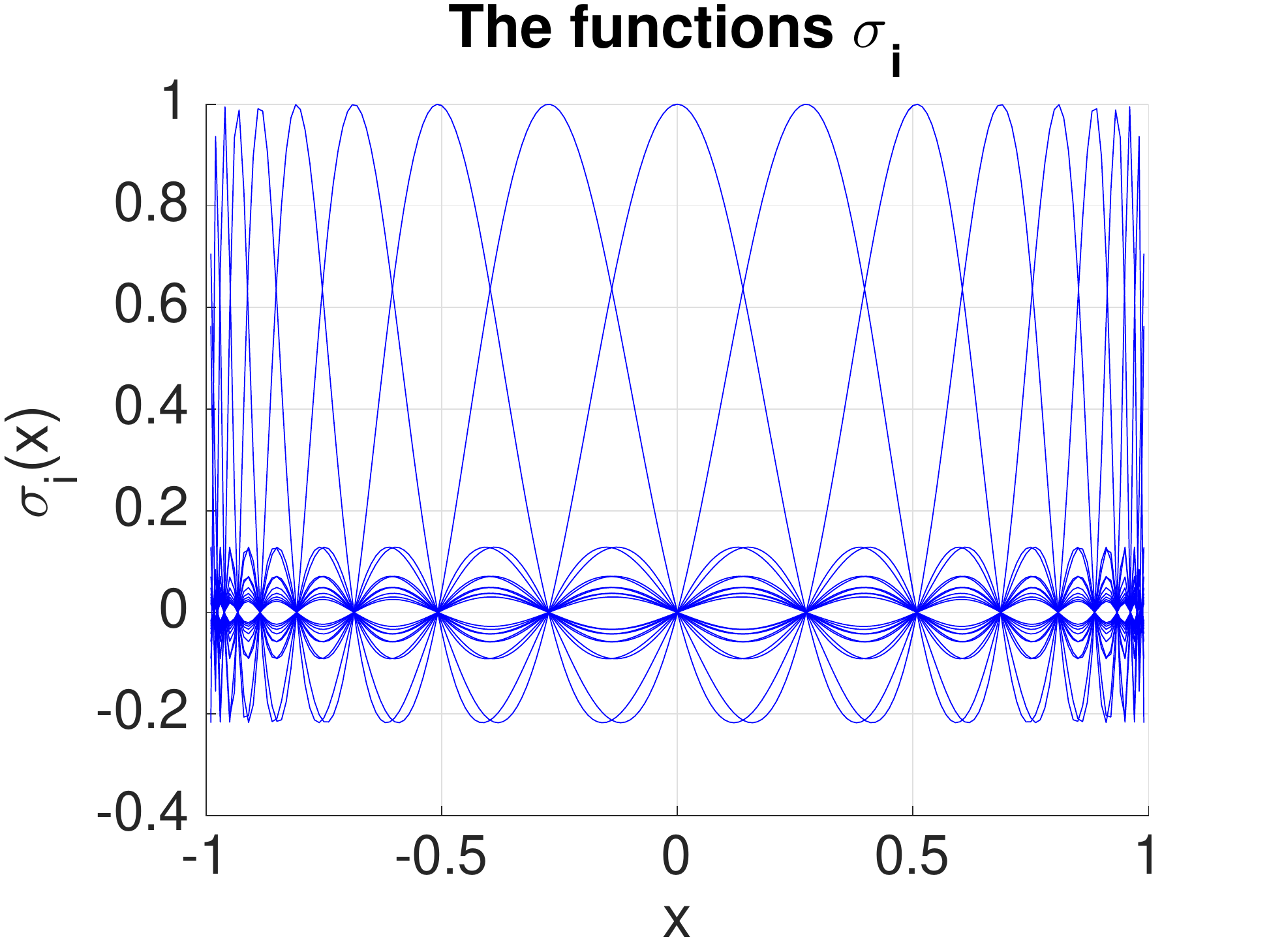}
	\caption{Extremal entries of the inverse collocation matrices $P_n^{-1}$ for the $\sinc$-basis from Sect.~\ref{exsinc} depending on $n$ (left) and the associate functions $\sigma_i:[-1,1]\to\R$, $-n\leq i\leq n$, for $n=10$ and $\tfrac{\delta}{\alpha}=1$ (right)}
	\label{figsinc}
\end{figure}
Collocation methods based on the \emph{cardinal sine function}
\begin{align*}
	\sinc:\R&\to\R,&
	\sinc x&:=
	\begin{cases}
		\tfrac{\sin(\pi x)}{\pi x},&x\neq 0,\\
		1,&x=0
	\end{cases}
\end{align*}
are discussed in \cite{stenger:93}. They require the $\tanh$-transformation $\phi:\R\to(a,b)$ 
\begin{align*}
	\phi(x)&:=\tfrac{a+b}{2}+\tfrac{b-a}{2}\tanh\intoo{\tfrac{x}{2}},&
	\phi^{-1}(x)&:=\ln\tfrac{x-a}{b-x}
\end{align*}
(with inverse). In order to approximate functions $u:[a,b]\to\R$ choose some $\delta\in(0,\pi)$ so that $(a,b)\subseteq D:=\bigl\{z\in\C:\,\bigl|\arg\tfrac{z-a}{b-z}\bigr|<\delta\bigr\}$ holds and assume $u$ is the restriction of a homomorphic function $u:D\to\C$ satisfying the growth condition $\abs{u(z)}\leq C\abs{(z-a)(z-b)}^\alpha$ for all $z\in D$ with reals $C,\alpha>0$. Given $n\in\N$, $h:=\sqrt{\tfrac{\pi\delta}{\alpha n}}$ and $S(j,h)(x):=\sinc\tfrac{x-jh}{h}$, we introduce the $d_n:=2n+3$ basis functions 
\begin{align*}
	\phi_{-(n+1)}(x)&:=\tfrac{b-x}{b-a},&
	\phi_i(x)&:=S(i,h)(\phi^{-1}(x))\fall\abs{i}\leq n,&
	\phi_{n+1}(x)&:=\tfrac{x-a}{b-a}
\end{align*}
with the $2n+3$ collocation points
$$
	x_j:=
	\begin{cases}
		a,&j=-n-1,\\
		\phi(jh),&\abs{j}\leq n,\\
		b,&j=n+1. 
	\end{cases}
$$
This yields the nonnegative collocation matrix $P_n:=(\phi_i(x_j))_{i,j=-n-1}^{n+1}$ given by
\begin{align*}
	P_n
	&=
	\begin{pmatrix}
		1 & \tfrac{b-x_{-n}}{b-a} & \ldots & \tfrac{b-x_n}{b-a} & 0\\
		0 & 1 & \ldots & 0 & 0\\
		\vdots & & \ddots & & \vdots\\
		0 & 0 & \ldots & 1 & 0\\
		0 & \tfrac{x_{-n}-a}{b-a} & \ldots & \tfrac{x_n-a}{b-a} & 1\\
	\end{pmatrix}\in\R^{(2n+3)\tm(2n+3)}.
\end{align*}
\fref{figsinc} shows the extremal entries of $P_n^{-1}$ depending on $n\geq 2$. Having positive and negative entries in $P_n^{-1}$, \tref{thmcol} fails to apply for $\sinc$-collocation.
\subsection{Bubnov-Galerkin methods}
For a second class of projection methods, assume $X(\Omega)=L^2(\Omega)$ and that $\set{\phi_1,\ldots,\phi_{d_n}}$ forms a basis of a subspace $X_n\subset L^2(\Omega)$. With the inner product
\begin{equation}
	(u,v):=\int_\Omega u(y)v(y)\d\mu(y)
	\label{noinner}
\end{equation}
we require orthogonality $(\pi_nu-u,\phi_i)=0$ for $1\leq i\leq d_n$ leading to functionals $\psi_i'(u)=(u,\phi_i)$. This results in the positively definite \emph{Gramian matrix} 
$$
	P_n=((\phi_j,\phi_i))_{i,j=1}^{d_n}. 
$$
\begin{thm}[positivity of $\Pi_n$ on $L^2(\Omega)^d$]
	Assume that every basis element $\phi_i:\Omega\to\R$, $1\leq i\leq d_n$, has nonnegative values. If the Gramian matrix $P_n$ is monomial, then $\Pi_n$ is $L^2(\Omega)_+^d$-positive. 
\end{thm}
\begin{pf}
	Let $u\in L^2(\Omega)^d$, $0\lle u$ and $y'\in Y_+'$. Thanks to \lref{lemconelp}(a) this implies $0\leq\sprod{u(y),y'}$ for $\mu$-a.a.\ $y\in\Omega$. Since the Gramian matrix $P_n$ is assumed to be monomial, we obtain from \cite[p.~2, Thm.~1.1]{kaczorek:02} that the inverse is nonnegative. Hence, if we write $\phi_j'(u)=((u_k,\phi_j))_{k=1}^d$, then \eqref{noproj} implies
	\begin{align*}
		\sprod{(\Pi_nu)(x),y'}
		&=
		\sprod{\sum_{j=1}^{d_n}\phi_j'(u)\phi_j(x),y'}
		=
		\sum_{j=1}^{d_n}\phi_j(x)\sprod{\phi_j'(u),y'}\\
		&=
		\sum_{j=1}^{d_n}\phi_j(x)\sprod{\sum_{i=1}^{d_n}(P_n^{-1})_{ij}\int_\Omega u(y)\phi_i(y)\d\mu(y),y'}\\
		&=
		\sum_{i=1}^{d_n}\int_\Omega \phi_i(y)\underbrace{\sprod{u(y),y'}}_{\geq 0}\d\mu(y)\sum_{j=1}^{d_n}\phi_j(x)(P_n^{-1})_{ij}\geq 0
	\end{align*}
	for all $x\in\Omega$ by assumption. Therefore, \lref{lemconelp}(a) yields the claim. 
	\qed
\end{pf}
\begin{cor}\label{corbg}
	For orthogonal bases $\set{\phi_1,\ldots,\phi_{d_n}}$ consisting of nonnegative functions the projection $\Pi_n$ is $L^2(\Omega)_+^d$-positive. 
\end{cor}
\begin{pf}
	Because $\set{\phi_1,\ldots,\phi_{d_n}}$ is orthogonal, the Gramian matrix is diagonal with positive entries along the diagonal, i.e.\ it is monomial. 
	\qed
\end{pf}
\begin{cor}
	Let Hypothesis $(L^2)$ hold with a monomial Gramian matrix $P_n$. If $\eK$ is $L^2(\Omega)_+^d$-positive, then also the compositions $\Pi_n\eK\in L(L^2(\Omega),X_n)$ and $\eK\Pi_n\in L(L^2(\Omega)^d$) are $L^2(\Omega)_+^d$-positive. 
\end{cor}
\begin{pf}
	We refer to \cref{corA2}. 
	\qed
\end{pf}

In the following subsections we endow $\Omega=[a,b]$ with a grid \eqref{grid}. 
\subsubsection{Piecewise constant functions in $\R^\kappa$}
\label{seccon}
The grid \eqref{grid} breaks $\Omega=[a,b]$ into the subintervals $[x_{j-1},x_j]$ of length $h_j$, $1\leq j\leq n$. For the ansatz space 
$$
	X_n:=\{u\in L^2[a,b]:\,u\text{ is constant on }[x_j,x_{j+1}], 0\leq j<n\}
$$
of piecewise constant functions the characteristic functions $\phi_j:=\chi_{[x_{j-1},x_j]}$, $1\leq j\leq n$, establish an orthogonal basis yielding the diagonal Gramian matrix 
$
	P_n:=\diag(h_0,\ldots,h_{n-1}).
$
In particular, \cref{corbg} applies and yields a positive Bubnov-Galerkin method with orthogonal projections
\begin{equation}
	\Pi_nu
	=
	\sum_{j=1}^n\frac{1}{h_{j-1}}\int_{x_{j-1}}^{x_j}u(y)\d y\chi_{[x_{j-1},x_j]}. 
	\label{pcproj}
\end{equation}
Concerning an extension to rectangles $\Omega:=[a_1,b_1]\tm\ldots\tm[a_\kappa,b_\kappa]$ we proceed as follows: Subdivide each interval $[a_j,b_j]$ into $n$ subintervals $I_i^j$. With the corresponding characteristic functions $\phi_i^j=\chi_{I_i^j}:[a_j,b_j]\to\R_+$, $1\leq i\leq n$, multivariable versions read as
$$
	\phi_\iota(x):=\prod_{j=1}^\kappa \phi_{\iota_j}^j(x_j)
	\text{ for }x=(x_1,\ldots,x_\kappa)\in\Omega,\,\iota=(\iota_1,\ldots,\iota_\kappa)\in\set{1,\ldots,n}^\kappa. 
$$
Given this, choose $\set{\phi_\iota:\Omega\to\R_+\mid\iota\in\set{1,\ldots,n}^\kappa}$ as basis of $X_n\subset L^p(\Omega)$ having the dimension $d_n=n^\kappa$. 

Differing from their convenient role in collocation methods (cf.\ \eref{expw}), piecewise linear functions do not preserve positivity.
\subsubsection{Piecewise linear functions}
\label{secplf}
\begin{figure}
	\includegraphics[scale=0.35]{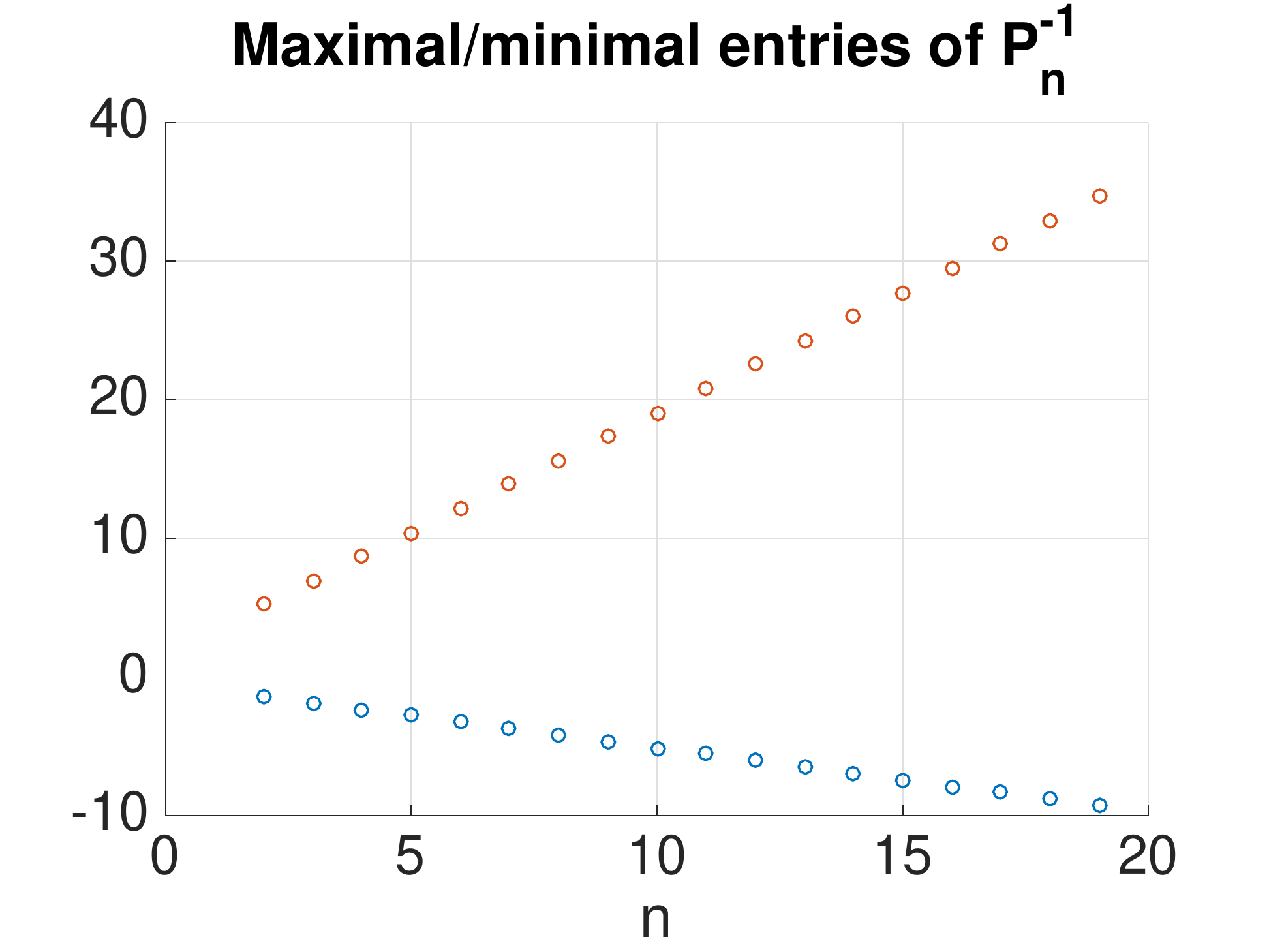}
	\includegraphics[scale=0.35]{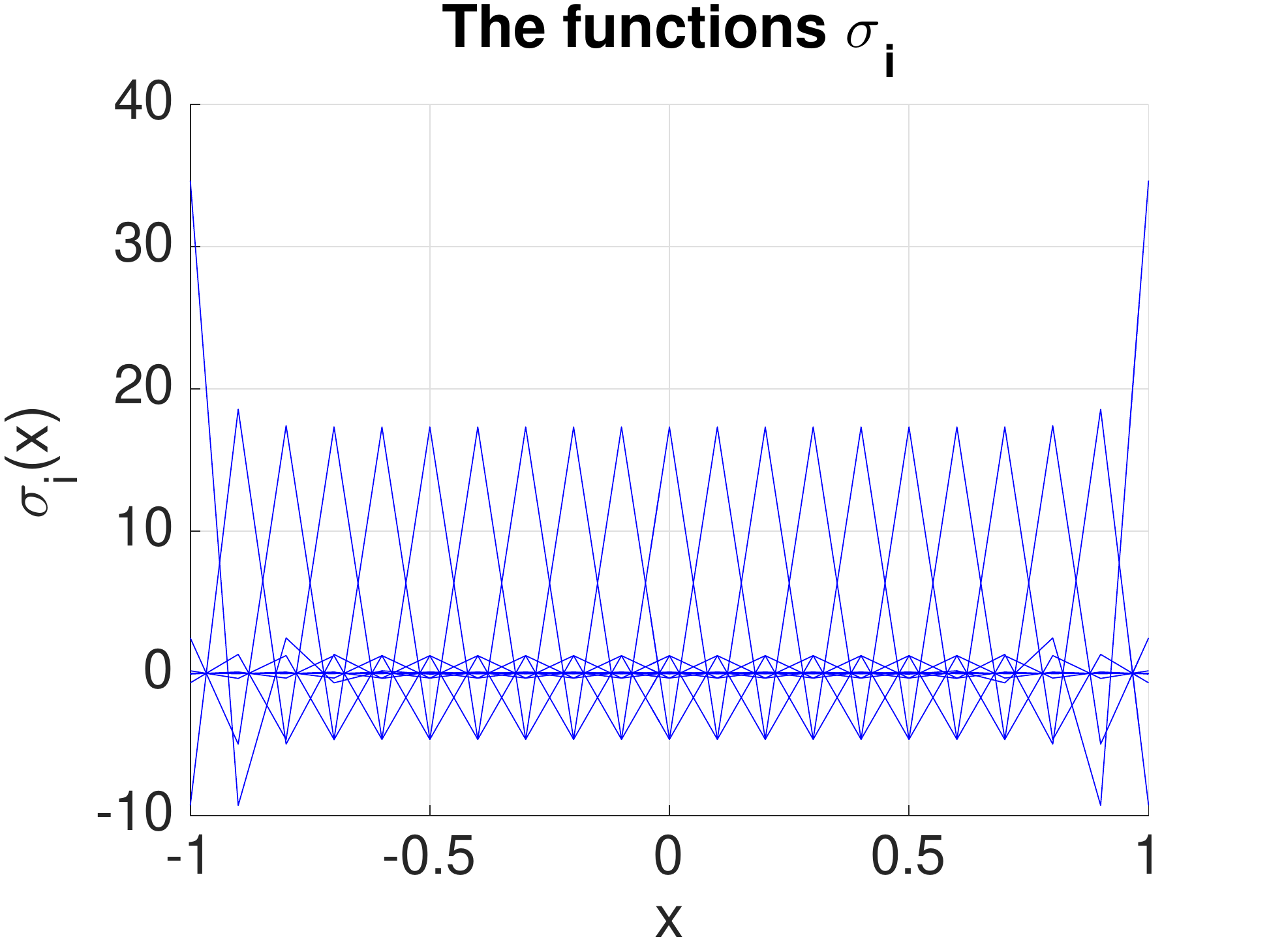}
	\caption{Extremal entries of the inverse Gramian matrices $P_n^{-1}$ for piecewise affine basis functions from Sect.~\ref{secplf} depending on $n$ (left), as well as the associate functions $\sigma_i:[-1,1]\to\R$, $0\leq i\leq n$, for $n=20$}
	\label{figlinear}
\end{figure}
With the above notation, consider the ansatz space 
$$
	X_n:=\set{u\in L^2[a,b]:\,u\text{ is affine on }[x_j,x_{j+1}], 0\leq j<n}
$$
of piecewise affine functions. The hat functions $\phi_j:[a,b]\to\R$ from \eref{expw} constitute a basis yielding the tridiagonal $(n+1)\tm(n+1)$-Gramian matrix 
$$
	P_n
	=
	\frac{1}{6}
	\begin{pmatrix}
		2h_0 & h_0 & & &\\
		h_0 & 2(h_0+h_1) & h_1 & &\\
		& \ddots & \ddots & \ddots & \\
		&& h_{n-2} & 2(h_{n-2}+h_{n-1}) & h_{n-1}\\
		&&& h_{n-1} & 2h_{n-1}
	\end{pmatrix}
$$
(cf.~\cite[p.~104, Rem.~4.5.26]{hackbusch:95}). Now \cref{corbg} does not apply, because as illustrated in \fref{figlinear}, the entries of $P_n^{-1}$ are positive and negative.
\begin{exmp}
	For simplicity we restrict to $h_j=\tfrac{b-a}{n}$ and obtain 
	\begin{equation*}
		P_n=
		\frac{b-a}{6n}
		\begin{pmatrix}
			2 & 1 & & & &\\
			1 & 4 & 1 & & &\\
			& \ddots & \ddots & \ddots & \\
			&& 1 & 4 & 1 \\
			&&& 1 & 2
		\end{pmatrix}.
	\end{equation*}
	as the Gramian matrix. For $x\in(a,a+h]$ one has
	\begin{align*}
		(\pi_nu)(x)
		\stackrel{\eqref{noproj}}{=} &
		\sum_{j=0}^n\phi_j'(u)\phi_j(x)
		=
		\phi_0'(u)\phi_0(x)+\phi_1'(u)\phi_1(x)\\
		\stackrel{\eqref{nolgs}}{=} &
		\intoo{(P_n^{-1})_{11}(u,\phi_0)+(P_n^{-1})_{12}(u,\phi_1)}\phi_0(x)\\
		&+
		\intoo{(P_n^{-1})_{21}(u,\phi_0)+(P_n^{-1})_{22}(u,\phi_1)}\phi_1(x)\\
		=&
		\intoo{(P_n^{-1})_{11}\phi_0(x)+(P_n^{-1})_{21}\phi_1(x)}(u,\phi_0)\\
		&+
		\intoo{(P_n^{-1})_{12}\phi_0(x)+(P_n^{-1})_{22}\phi_1(x)}(u,\phi_1)
	\end{align*}
	and in particular with the characteristic function $u:=\chi_{[a,a+\theta h]}\in L^p[a,b]_+$ and $\theta\in(0,1)$ it results from \eqref{noinner} that $(u,\phi_0)=\tfrac{h\theta}{2}(2-\theta)$, $(u,\phi_1)=\tfrac{h\theta}{2}\theta$. We consequently arrive at
	\begin{align*}
		(\pi_nu)(x)
		=&
		\tfrac{h\theta}{2}\left[
		\intoo{(P_n^{-1})_{11}\phi_0(x)+(P_n^{-1})_{21}\phi_1(x)}(2-\theta)
		\right.\\
		&+
		\left.
		\intoo{(P_n^{-1})_{12}\phi_0(x)+(P_n^{-1})_{22}\phi_1(x)}\theta
		\right]
	\end{align*}
	and hence, if we choose $\theta>0$ close to $0$ and $x\leq a+h$ close to $a+h$ (where $\phi_0$ has values near $0$ and $\phi_1$ near 1, cf.~\eqref{hatdef}), then the sign of $(\pi_nu)(x)$ is determined by the entry $(P_n^{-1})_{21}$. Now by \lref{leminv} one has $(P_n^{-1})_{21}<0$ and thus there exists an interval $[a+h_0,a+h]$, $h_0\in(0,h)$ on which $(\pi_n u)(x)$ is negative, that is $\pi_nu\not\in L^2[a,b]_+$. 
\end{exmp}
\section{Applications}
\label{sec5}
\subsection{Dispersal kernels}
In order to describe the dispersal stage in ecological models \cite{lutscher:19}, various real-valued functions $k_{ij}:\Omega^2\to\R_+$ are used as entries in corresponding matrix-valued kernels $K:\Omega^2\to\R^{d\tm d}$ for Fredholm operators $\eK$; thus one speaks of \emph{dispersal kernels} and chooses the Lebesgue measure $\mu=\lambda_\kappa$. Often the functions $k=k_{ij}$ are symmetric, i.e.\ $k(x,y)=k(y,x)$, or of convolution form $k(x,y)=\tilde k(x-y)$ for all $x,y\in\Omega$ with an even function $\tilde k:\R^\kappa\to\R_+$. Given a \emph{dispersal rate} $\alpha>0$ and the Euclidean norm $\abs{\cdot}$ on $\R^\kappa$ typical choices are listed in Tab.~\ref{tabkernel}. 
\begin{table}[ht]
	\begin{tabular}{r|l}
		kernel & $\tilde k(x)$\\ 
		\hline
		Gau{\ss} & $\tfrac{1}{\sqrt{2\pi\alpha^2}}\exp\intoo{-\tfrac{1}{2\alpha^2}\abs{x}^2}$\\ 
		Cauchy & $\tfrac{\alpha}{\pi(\alpha^2+\abs{x}^2)}$\\ 
		Laplace & $\tfrac{1}{2\alpha}\exp\intoo{-\tfrac{1}{\alpha}\abs{x}}$\\ 
		exponential square root & $\frac{1}{4\alpha}\exp\intoo{-\sqrt{\tfrac{1}{\alpha}\abs{x}}}$\\
		top-hat & $\tfrac{1}{2\alpha}\chi_{B_{\alpha}(0)}(\abs{x})$\\ 
		tent & $\tfrac{1}{\alpha}\max\set{0,1-\tfrac{1}{\alpha}\abs{x}}$\\ 
	\end{tabular}
	\caption{List of dispersal kernels satisfying $\int_{\R}\tilde k(y)\d y=1$ for $\kappa=1$}
	\label{tabkernel}
\end{table}
The kernels from Tab.~\ref{tabkernel} satisfy Hypothesis $(L)$ (whose limit relation moreover holds uniformly in $x_0$), as well as Hypothesis $(L^2)$. Consequently, the resulting Fredholm operator $\eK u=\int_\Omega k(\cdot,y)u(y)\d y$ acting on $C(\Omega)$ resp.\ $L^2(\Omega)$ is even compact. In case $Y_+=\R_+$ these kernels yield positive operators $\eK$ due to \tref{thmlinmon} resp.\ \ref{thml2}. On the space $X(\Omega)=C(\Omega)$, besides for the top-hat and the tent kernel, $\eK$ is even strongly positive. Yet, for sufficiently large dispersal rates $\alpha$ also the top-hat and tent kernel yield a strongly positive $\eK$, i.e.\ provided the condition $\Omega\subseteq B_\alpha(y)$ for all $y\in\Omega$ holds. 

Now we return to matrix-valued kernels in $\eK$: Here criteria for positivity can be derived from \eref{exgencone} and \ref{exorthants}. With vectors $e_i\in\R^d$ spanning $Y_+\subset\R^d$ and $e_j'\in\R^d$, a kernel $K(x,y)\in\R^{d\tm d}$ is
\begin{itemize}
	\item $Y_+$-positive $\Leftrightarrow 0\leq\sprod{K(x,y)e_i,e_j'}$, 

	\item strictly $Y_+$-positive $\Leftrightarrow 0\leq\sprod{K(x,y)e_i,e_j'}$ and $K(x,y)$ is $Y_+$-injective,

	\item strongly $Y_+$-positive $\Leftrightarrow 0<\sprod{K(x,y)e_i,e_j'}$ 
\end{itemize}
holds for all $1\leq i,j\leq d$, $x\in\Omega$ and $\mu$-a.a.\ $y\in\Omega$. 

If a cone $Y_+$ is an orthant as in \eref{exorthants}, then $\sprod{K(x,y)e_i,e_j'}=\varsigma_i\varsigma_jk_{ij}(x,y)$ implies that $K(x,y)$ is
\begin{itemize}
	\item $Y_+$-positive $\Leftrightarrow\varsigma_i\varsigma_jk_{ij}(x,y)\geq 0$, 

	\item strictly $Y_+$-positive $\Leftrightarrow\varsigma_i\varsigma_jk_{ij}(x,y)\geq 0$ and $K(x,y)$ is $Y_+$-injective, 

	\item strongly $Y_+$-positive $\Leftrightarrow\sgn k_{ij}(x,y)=\varsigma_i\varsigma_j$ 
\end{itemize}
for all $1\leq i,j\leq d$, $x\in\Omega$ and $\mu$-a.a.\ $y\in\Omega$. 
\subsection{Krein-Rutman eigenfunctions}
Let $X(\Omega)$ be a Banach space of real-valued functions over $\Omega$. Given a compact Fredholm operator $\eK\in L(X(\Omega)^d)$, the Krein-Rutman \tref{thmkr} provides conditions guaranteeing that its spectral radius $r(\eK)>0$ is the dominant eigenvalue with an associate nonnegative eigenfunction $u^\ast\in X(\Omega)^d$. In this section, we investigate in which sense this property is preserved under discretization. 

For this purpose, consider an abstract eigenvalue problem
\begin{equation}
	\eK u=\lambda u\quad\text{in }X(\Omega)^d
	\label{evp}
\end{equation}
for a Fredholm operator $\eK$. Numerical approaches to \eqref{evp} approximate $(\lambda,u)\in\C\tm X(\Omega)$ by eigenpairs $(\lambda_n,\hat v^n)\in\C\tm\R^{dN_n}$ of related problems
\begin{equation}
	K^n\hat v^n=\lambda_n\hat v^n,
	\label{evpd}
\end{equation}
where the matrices $K^n\in\R^{dN_n\tm dN_n}$ depend on the particular discretization method. Nevertheless, $K^n$ consists of block matrices
\begin{align*}
	K^n
	&:=
	\begin{pmatrix}
		K_{11} & \ldots & K_{1d}\\
		\vdots & & \vdots\\
		K_{d1} & \ldots & K_{dd}
	\end{pmatrix},&
	\hat v^n
	&:=
	\begin{pmatrix}
		v_1\\
		\vdots\\
		v_d
	\end{pmatrix},
\end{align*}
with also $v_i\in\R^{N_n}$ and $K_{ij}\in\R^{N_n\tm N_n}$ depending on the method. 

Nevertheless, in our numerical computations we rely on the \textsc{Matlab} functions \texttt{eig} to approximate the dominant eigenvalue of $K^n$ and \texttt{eigs} when also the corresponding eigenvectors are of interest. 
\subsubsection{Nystr\"om methods}
Let $X(\Omega)=C(\Omega)$ and \eqref{quad} be a quadrature rule. Replacing $\eK$ in \eqref{evp} by the discrete Fredholm operator \eqref{nofredn} yields
$$
	\sum_{j_2=0}^{q_n}w_{j_2}K(x,\eta_{j_2})u(\eta_{j_2})=\lambda u(x)\fall x\in\Omega. 
$$
If we now set $x=\eta_{j_1}$ for $0\leq j_1\leq q_n$, then one obtains a matrix $K^n$ as in relation \eqref{evpd} with $N_n:=q_n+1$ and the blocks
$$
	K_{i_1i_2}=(w_{j_2}k_{i_1i_2}(\eta_{j_1},\eta_{j_2}))_{j_1,j_2=0}^{q_n}
	\fall 1\leq i_1,i_2\leq d. 
$$
For simplicity we retreat to real-valued continuous kernels $K:\Omega\tm\Omega\to(0,\infty)$ and the solid cone $Y_+:=\R_+$. Then $\eK$ is strongly positive by \tref{thmlinmon}(b) and the Krein-Rutman \tref{thmkr}(b) yields that the spectral radius $r(\eK)>0$ is the dominant eigenvalue with eigenfunction $u^\ast\in C(\Omega)_+^\circ$. If all the quadrature weights are positive, then a discrete Fredholm operator $\eK^n$ from \eqref{nofredn} is
\begin{itemize}
	\item strongly positive on $C(\Omega_n)$ (cf.~\rref{remsame}) and thus \tref{thmkr}(b) applies, 

	\item positive on $C(\Omega)$ (due to \tref{thmkny}), where \tref{thmkr}(a) applies.
\end{itemize}
In any case, one obtains the dominant eigenpair $(r(\eK^n),\hat v^n)$ from the equation \eqref{evpd} and as dominant eigenfunction $u_n^\ast\in C(\Omega)$ we use the \emph{Nystr\"om interpolate}
$$
	u_n^\ast(x):=\frac{1}{r(\eK^n)}\sum_{j=0}^{q_n}w_j K(x,\eta_j)\hat v_j^n\fall x\in\Omega.
$$

\begin{figure}
	\includegraphics[scale=0.35]{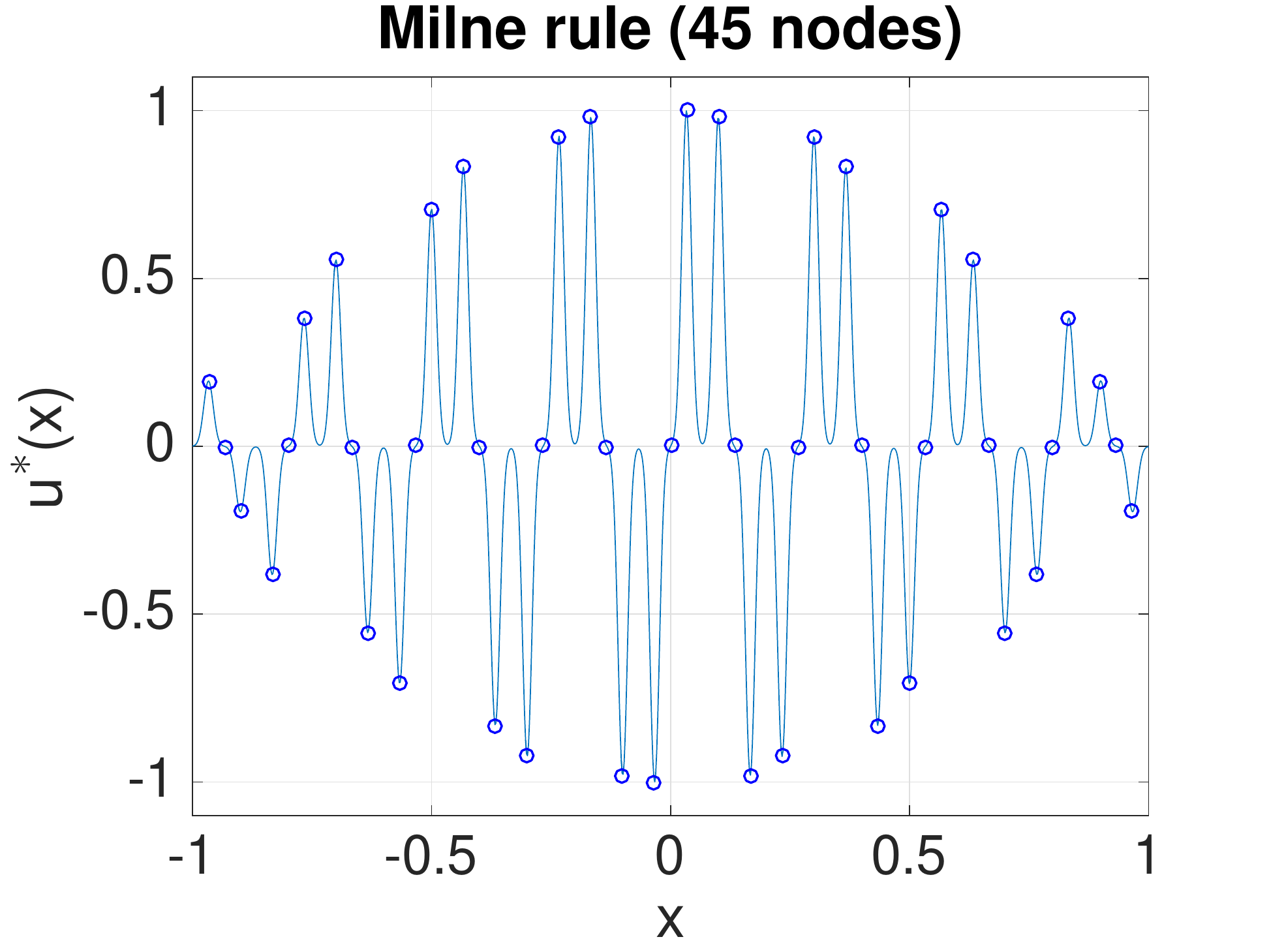}
	\includegraphics[scale=0.35]{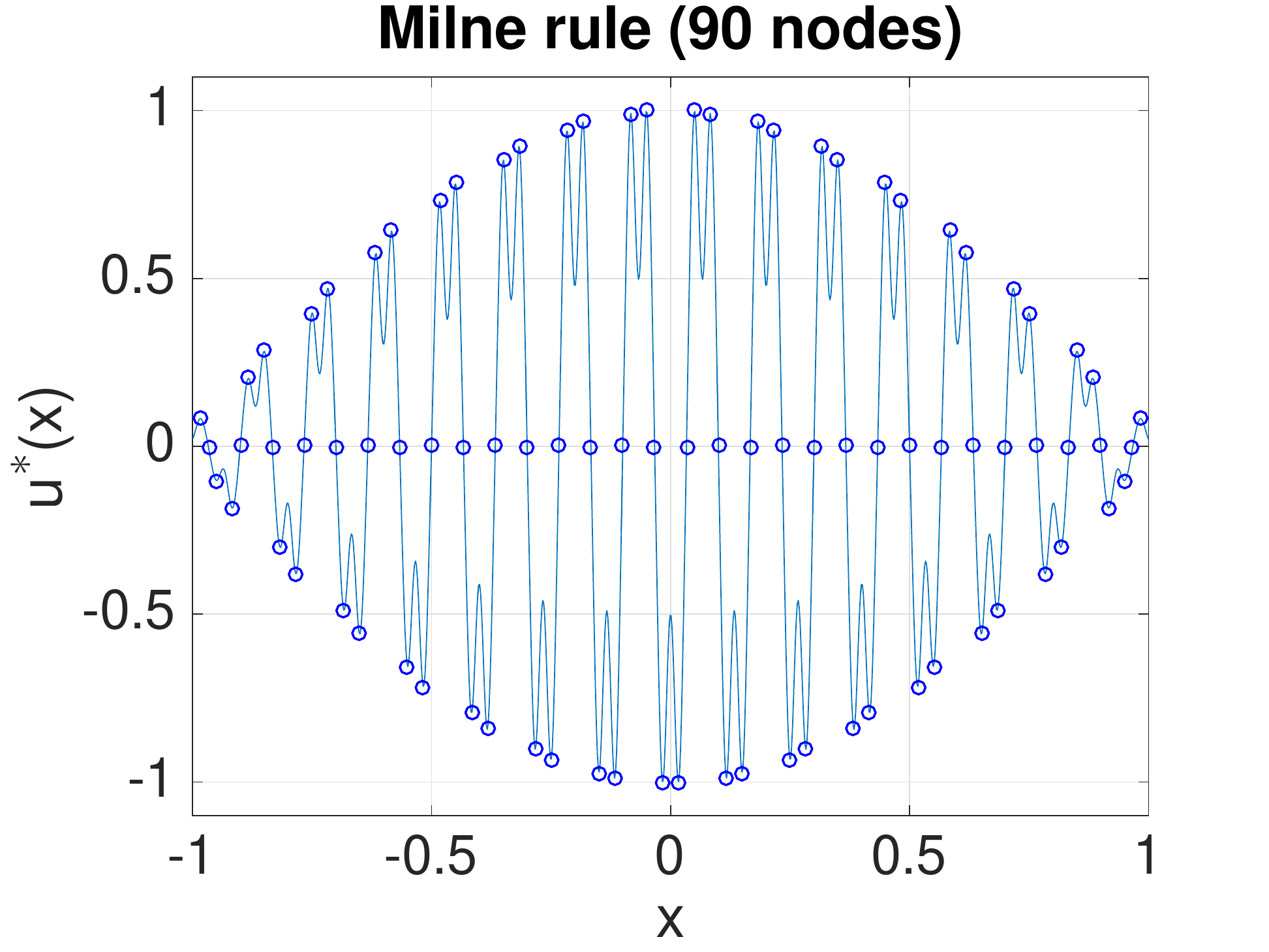}\\
	\caption{Nystr\"om discretization to approximate the dominant positive eigenfunction of the Gau{\ss} kernel with $\alpha=0.01$ based on Milne's rule \eqref{nomilne} with with $45$ nodes (left) and $90$ nodes (right).}
	\label{figkr2}
\end{figure}

The positivity of the numerically obtained eigenfunction $u_n^\ast$ is confirmed by 
\begin{figure}
	\includegraphics[scale=0.35]{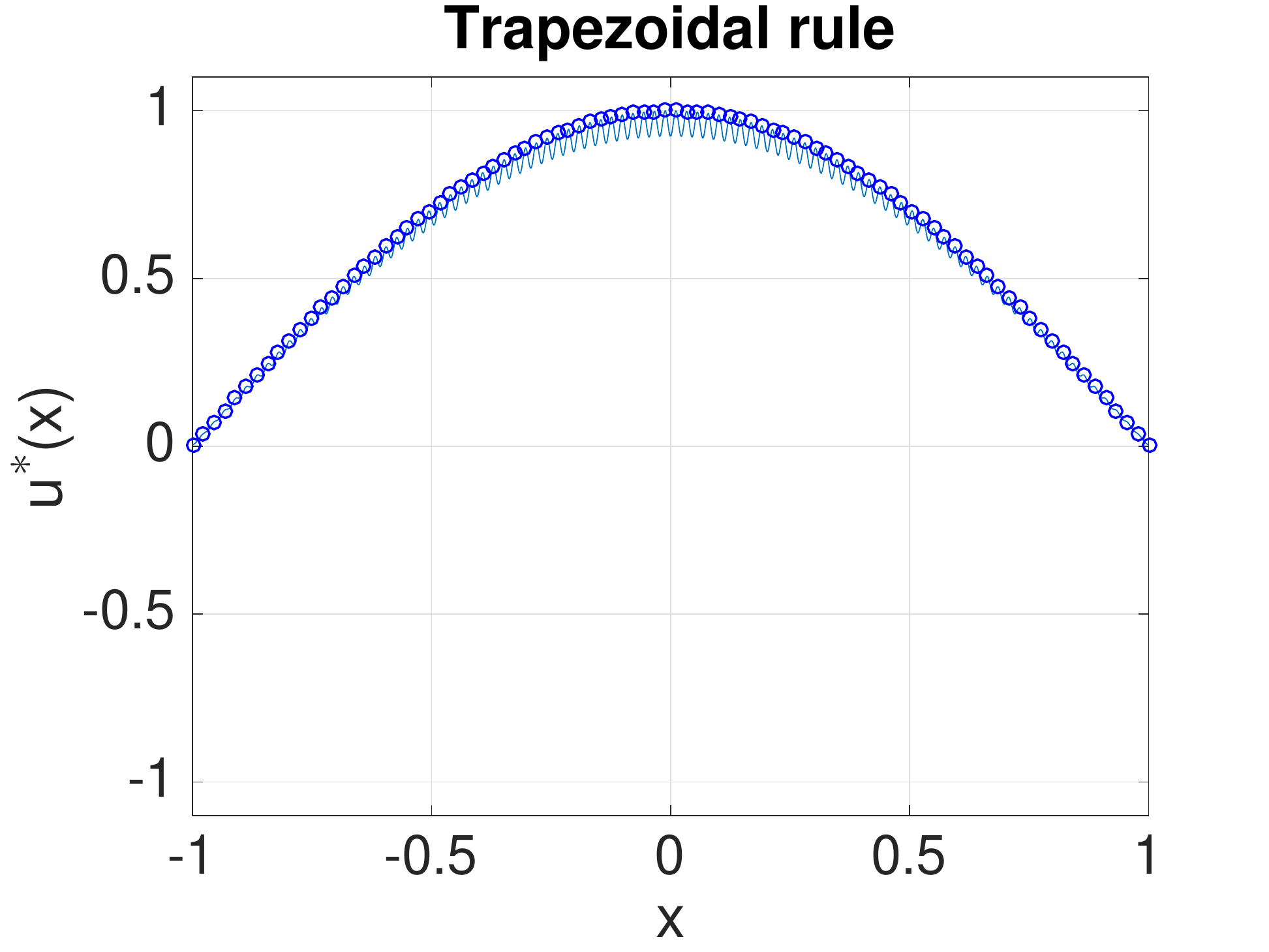}
	\includegraphics[scale=0.35]{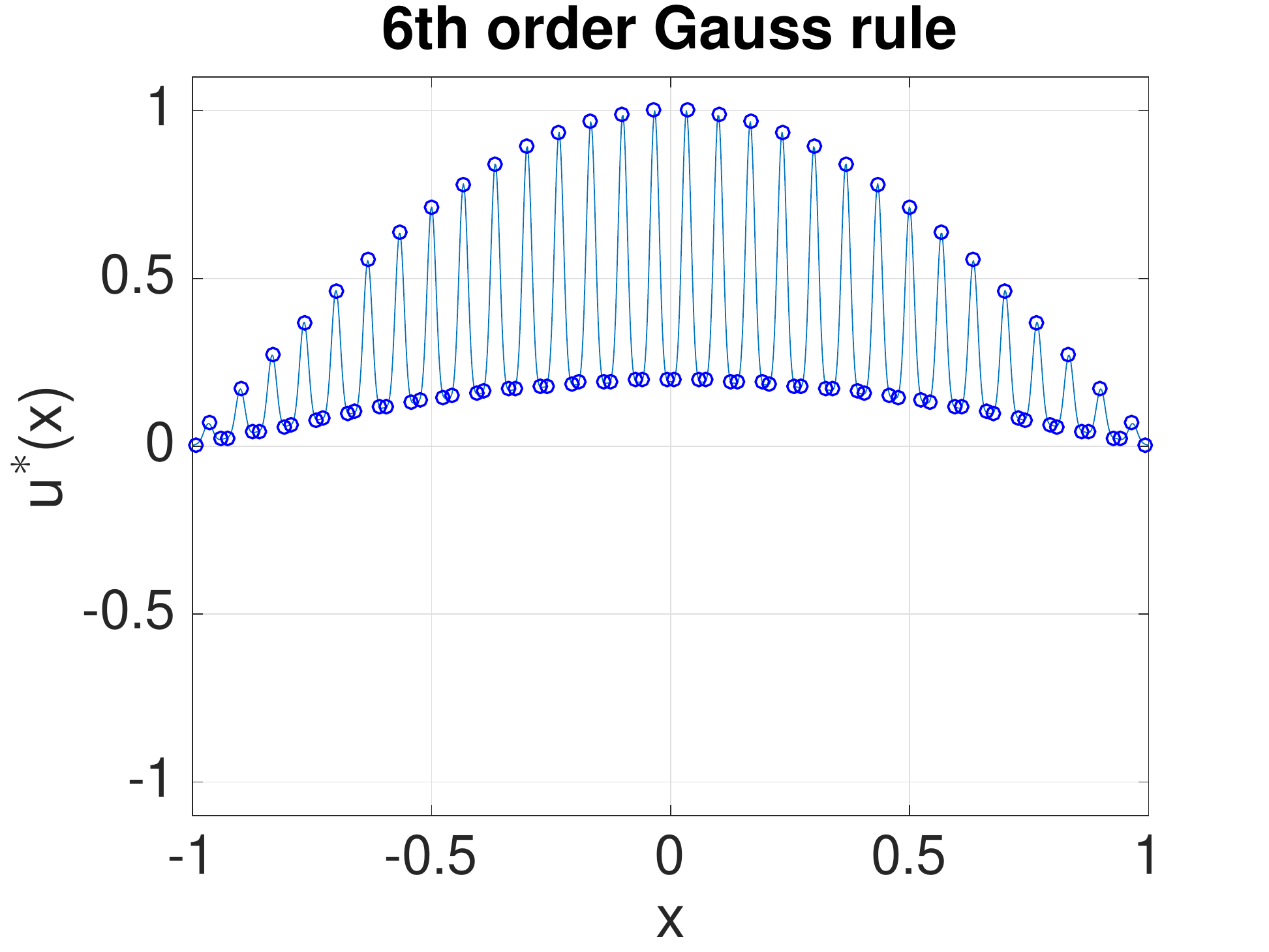}\\
	\caption{Nystr\"om discretization to approximate the dominant positive eigenfunction of the Gau{\ss} kernel with $\alpha=0.01$ based on quadrature rules with $90$ nodes. The dots represent the values $u_n^\ast(\eta)$, $\eta\in\Omega_n$}
	\label{figkr1}
\end{figure}
\begin{exmp}[Gau{\ss} kernel]
	We apply different Nystr\"om methods to the Fredholm operator
	$$
		(\eK u)(x)
		=
		\frac{1}{\sqrt{2\pi\alpha^2}}
		\int_{-1}^1e^{-\tfrac{(x-y)^2}{2\alpha^2}}u(y)\d y		
	$$
	with the dispersal rate $\alpha=0.01$ over $\Omega=[-1,1]$. The corresponding dominant eigenfunctions (as Nystr\"om interpolates) for the Trapezoidal rule \eqref{notrap} and $6$th order Gau{\ss} rule \eqref{nogauss} are illustrated in \fref{figkr1} (each with $90$ notes). Both cases yield strongly positive eigenfunctions (although hardly visible in \fref{figkr1}, one has $u^\ast(x)>0$ for $x\in\set{-1,1}$), but the Gau{\ss} discretizations exhibit a strongly oscillatory behavior. This is due to the fact that the used number of nodes $(q_n=89)$ is small in comparison to the dispersal rate $(\alpha=0.01)$. Oscillations become weaker when increasing $q_n$.
	\\
	For the Milne rule \eqref{nomilne} this situation changes, since it involves negative weights. Indeed, the resulting discrete Fredholm operator $\eK^n$ is not positive and the Krein-Rutman \tref{thmkr} does not apply. As illustrated in \fref{figkr2} the eigenfunction corresponding to the eigenvalue with maximal real part can have varying signs. Numerical experiments exhibit that this sign changing property persists until the number of nodes is larger than $146$. 
\end{exmp}
\subsubsection{Projection methods}
Based on the ansatz space $X_n=\spann\set{\phi_1,\ldots,\phi_{d_n}}$ we aim to approximate solutions to eigenvalue problems \eqref{evp} numerically. Thereto, replace $u$ in \eqref{evp} by $u^n=\Pi_nu$ given as in \eqref{no41s}. This leads to a discretized version
\begin{equation}
	\eK^n u^n=\lambda_n u^n\quad\text{in }X_n^d
	\label{evpn}
\end{equation}
with $\eK^n=\eK\Pi_n$ in order to obtain approximations $(\lambda_n,u^n)$, $u^n\in X_n^d$, to an eigenpair $(\lambda,u)$, $u\in X(\Omega)^d$, of \eqref{evp}. Note that for a solution $u^n\in X_n^d$ of the discretized problem \eqref{evpn} one has $\eK^nu^n\in X_n^d$. 

Similarly, one proceeds for spatial discretization of Fredholm equations of the first and second kind (cf.~\cite{atkinson:97,hackbusch:95}). 

For a matrix-valued kernel $K(x,y)$ define the Fredholm operators
\begin{align*}
	\eK_{i_1,i_2}&\in L(X(\Omega)),&
	\eK_{i_1i_2}u&:=\int_\Omega k_{i_1i_2}(\cdot,y)u(y)\d y\fall 1\leq i_1,i_2\leq d. 
\end{align*}
In this notation, the eigenvalue problem \eqref{evp} is equivalent to the relation
$
	\sum_{i_2=1}^d
	\eK_{i_1i_2}u_{i_2}
	=
	\lambda u_{i_1}
$
for $1\leq i_1\leq d$ and in order to arrive at the spatially discretized problem \eqref{evpn}, we insert $u_{i_2}=\sum_{j_1=1}^{d_n}v_{i_2}^{j_1}\phi_{j_1}$ with $v_{i_2}^{j_1}\in\R$. Applying the functional $\phi_{j_2}'$ yields
$$
	\sum_{i_2=1}^d\sum_{j_1=1}^{d_n}
	\phi_{j_2}'(\eK_{i_1i_2}\phi_{j_1})v_{i_2}^{j_1}
	=
	\lambda_n v_{i_1}^{j_2}
	\fall 1\leq i_1\leq d,\,1\leq j_2\leq d_n.
$$
This is equivalent to an eigenproblem \eqref{evpd} with $N_n=d_n$ and the blocks
\begin{align*}
	K_{i_1i_2}
	&:=
	\begin{pmatrix}
		\phi_1'(\eK_{i_1i_2}\phi_1) & \ldots & \phi_1'(\eK_{i_1i_2}\phi_{d_n})\\
		\vdots & & \vdots\\
		\phi_{d_n}'(\eK_{i_1i_2}\phi_1) & \ldots & \phi_{d_n}'(\eK_{i_1i_2}\phi_{d_n})
	\end{pmatrix},&
	v_{i_1}
	&:=
	\begin{pmatrix}
		v_{i_1}^1\\
		\vdots\\
		v_{i_1}^{d_n}
	\end{pmatrix}
\end{align*}
for all $1\leq i_1,i_2\leq d$. On this basis, the eigenfunctions of \eqref{evpn} become
$$
	u^n=
	\sum_{j=1}^{d_n}
	\phi_j
	\begin{pmatrix}
		v_1^j\\
		\vdots\\
		v_d^{j}
	\end{pmatrix}.
$$
In order to be more specific, let us emphasize the subsequent projection methods with an interval $\Omega=[a,b]$ on a grid \eqref{grid}. 

\emph{Collocation with Lagrange bases}: 
Lagrange bases $\set{\phi_1,\ldots,\phi_{d_n}}$ lead to blocks
$$
	K_{i_1i_2}
	=
	\int_\Omega
	\begin{pmatrix}
		k_{i_1i_2}(x_1,y)\phi_1(y) & \ldots & k_{i_1i_2}(x_1,y)\phi_{d_n}(y)\\
		\vdots & & \vdots\\
		k_{i_{1}i_2}(x_{d_n},y)\phi_1(y) & \ldots & k_{i_1i_2}(x_{d_n},y)\phi_{d_n}(y)\\
	\end{pmatrix}
	\d y.
$$
In particular, the bases from \eref{expw} (hat functions) and \eref{exlag} (Lagrange functions) are of the form $\set{\phi_0,\ldots,\phi_n}$ and thus $d_n=n+1$ in \eqref{evpd}. 

\emph{Collocation with cubic splines}: 
Use the order-preserving $C^1$-spline 
$$
	(\Pi_nu)(x)
	:=
	\kappa_j(x)u(x_j)+\bar\kappa_j(x)u(x_{j+1})
	\fall x\in[x_j,x_{j+1}],\,0\leq j<n
$$
from \cref{corspline} with the real coefficients
\begin{align*}
	\kappa_j(x)&:=1-3(\tfrac{x-x_j}{h_j})^2+2(\tfrac{x-x_j}{h_j})^3,&
	\bar\kappa_j(x)&:=3(\tfrac{x-x_j}{h_j})^2-2(\tfrac{x-x_j}{h_j})^3,
\end{align*}
which results in
\begin{align*}
	(\eK\Pi_nu)(x)
	&=
	\sum_{j=0}^{n-1}\int_{x_j}^{x_{j+1}}
	K(x,y)\kappa_j(y)\d yu(x_j)
	+
	\sum_{j=0}^{n-1}\int_{x_j}^{x_{j+1}}
	K(x,y)\bar\kappa_j(y)\d yu(x_{j+1})
\end{align*}
for all $x\in[a,b]$. 
Setting $x=x_i$ for $0\leq i\leq n$ in the ansatz \eqref{evpn} combined with the abbreviations $u_i:=u(x_i)$ yields the $n+1$ identities
\begin{align*}
	\sum_{j=0}^{n-1}\int_{x_j}^{x_{j+1}}K(x_i,y)\kappa_j(y)\d yu_j
	+
	\sum_{j=0}^{n-1}\int_{x_j}^{x_{j+1}}K(x_i,y)\bar\kappa_j(y)\d yu_{j+1}
	&=
	\lambda u_i
\end{align*}
for all $0\leq i\leq n$, which are equivalent to the eigenvalue problem
\begin{align*}
	&
	\Bigl(
	\begin{pmatrix}
		\int_{x_0}^{x_1}K(x_0,y)\kappa_0(y)\d y & 
		\ldots & 
		\int_{x_{n-1}}^{x_n}K(x_0,y)\kappa_{n-1}(y)\d y & 
		0\\
		\vdots & & \vdots & \vdots\\
		\int_{x_0}^{x_1}K(x_{n-1},y)\kappa_0(y)\d y & 
		\ldots & 
		\int_{x_{n-1}}^{x_n}K(x_{n-1},y)\kappa_{n-1}(y)\d y & 
		0\\
		0 & \ldots & 0 & 0
	\end{pmatrix}\\
	&+
	\begin{pmatrix}
		0 & 0 & \cdots & 0\\
		0 & \int_{x_0}^{x_1}K(x_0,y)\bar\kappa_0(y)\d y & 
		\ldots & 
		\int_{x_{n-1}}^{x_n}K(x_0,y)\bar\kappa_{n-1}(y)\d y\\
		\vdots & \vdots & & \vdots\\
		0 & \int_{x_0}^{x_1}K(x_{n-1},y)\bar\kappa_0(y)\d y & 
		\ldots & 
		\int_{x_{n-1}}^{x_n}K(x_{n-1},y)\bar\kappa_{n-1}(y)\d y\\
	\end{pmatrix}
	\Bigr)
	\begin{pmatrix}
		u_0\\
		\vdots\\
		u_{n-1}\\
		u_n
	\end{pmatrix}
	=
	\lambda
	\begin{pmatrix}
		u_0\\
		\vdots\\
		u_{n-1}\\
		u_n
	\end{pmatrix}
\end{align*}
in $\R^{(n+1)d}$. Referring to the formulation \eqref{evpd} this gives rise to the blocks
\begin{align*}
	K_{i_1i_2}
	=
	&
	\begin{pmatrix}
		\int_{x_0}^{x_1}k_{i_1i_2}(x_0,y)\kappa_0(y)\d y & 
		\ldots & 
		\int_{x_{n-1}}^{x_n}k_{i_1i_2}(x_0,y)\kappa_{n-1}(y)\d y & 
		0\\
		\vdots & & \vdots & \vdots\\
		\int_{x_0}^{x_1}k_{i_1i_2}(x_{n-1},y)\kappa_0(y)\d y & 
		\ldots & 
		\int_{x_{n-1}}^{x_n}k_{i_1i_2}(x_{n-1},y)\kappa_{n-1}(y)\d y & 
		0\\
		0 & \ldots & 0 & 0
	\end{pmatrix}\\
	&+
	\begin{pmatrix}
		0 & 0 & \cdots & 0\\
		0 & \int_{x_0}^{x_1}k_{i_1i_2}(x_0,y)\bar\kappa_0(y)\d y & 
		\ldots & 
		\int_{x_{n-1}}^{x_n}k_{i_1i_2}(x_0,y)\bar\kappa_{n-1}(y)\d y\\
		\vdots & \vdots & & \vdots\\
		0 & \int_{x_0}^{x_1}k_{i_1i_2}(x_{n-1},y)\bar\kappa_0(y)\d y & 
		\ldots & 
		\int_{x_{n-1}}^{x_n}k_{i_1i_2}(x_{n-1},y)\bar\kappa_{n-1}(y)\d y\\
	\end{pmatrix}. 
\end{align*}

\emph{Bubnov-Galerkin method with piecewise constant basis}:
One obtains the blocks
$$
	K_{i_1i_2}
	=
	\int_\Omega\int_\Omega
	\begin{pmatrix}
		k_{i_1i_2}(x,y)\phi_1(y)\phi_1(x) & \ldots & k_{i_1i_2}(x,y)\phi_{d_n}(y)\phi_1(x)\\
		\vdots & & \vdots\\
		k_{i_1i_2}(x,y)\phi_1(y)\phi_{d_n}(x) & \ldots & k_{i_1i_2}(x,y)\phi_{d_n}(y)\phi_{d_n}(x)
	\end{pmatrix}
	\d y
	\d x.
$$
Because the piecewise constant basis functions $\set{\phi_1,\ldots,\phi_n}$ from Sect.~\ref{seccon} have the support $[x_j,x_{j+1}]$, $0\leq j<n$, the integrals simplify to
$$
	K_{i_1i_2}
	=
	\begin{pmatrix}
		\int_{x_0}^{x_1}\int_{x_0}^{x_1} k_{i_1i_2}(x,y)\d y\d x & \ldots & \int_{x_0}^{x_1}\int_{x_{n-1}}^{x_n} k_{i_1i_2}(x,y)\d y\d x\\
		\vdots & & \vdots\\
		\int_{x_{n-1}}^{x_n}\int_{x_0}^{x_1} k_{i_1i_2}(x,y)\d y\d x & \ldots & \int_{x_{n-1}}^{x_n}\int_{x_{n-1}}^{x_n} k_{i_1i_2}(x,y)\d y\d x
	\end{pmatrix}.
$$
\begin{table}
	\begin{tabular}{r|lll}
		$i$ & $\alpha_i$ & $\nu_i$ & $\lambda^i$ \\
		\hline
		$1$ & 1.0 & 0.86033358901938 & 0.5746552163364324\\
		$2$ & 2.0 & 1.30654237418881 & 0.3694054047082261
	\end{tabular}
	\caption{Dominant eigenvalues $\lambda^i$ of $\eK_i$ for $L=2$}
	\label{tabev}
\end{table}

\begin{figure}
	\includegraphics[scale=0.35]{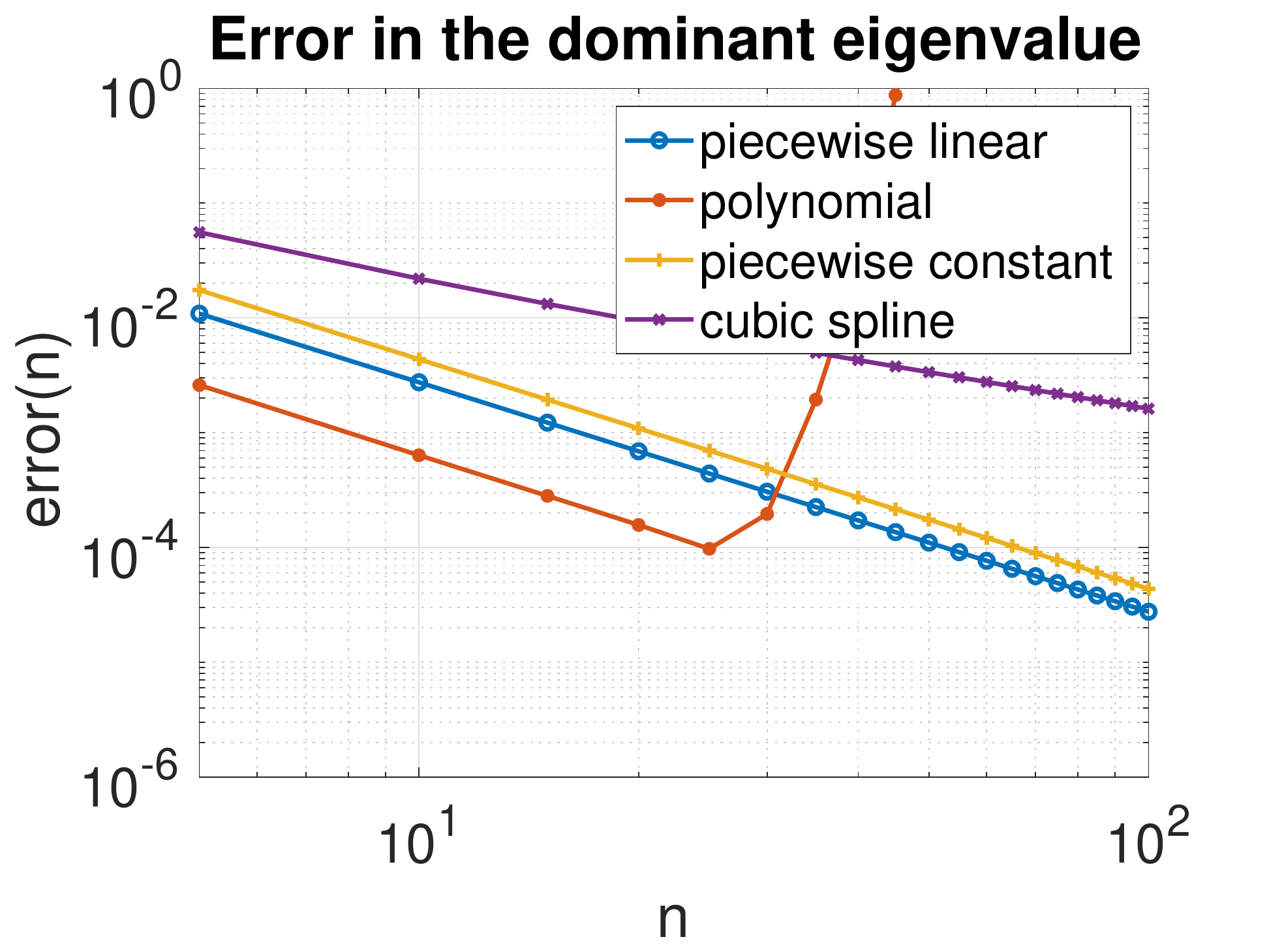}
	\includegraphics[scale=0.35]{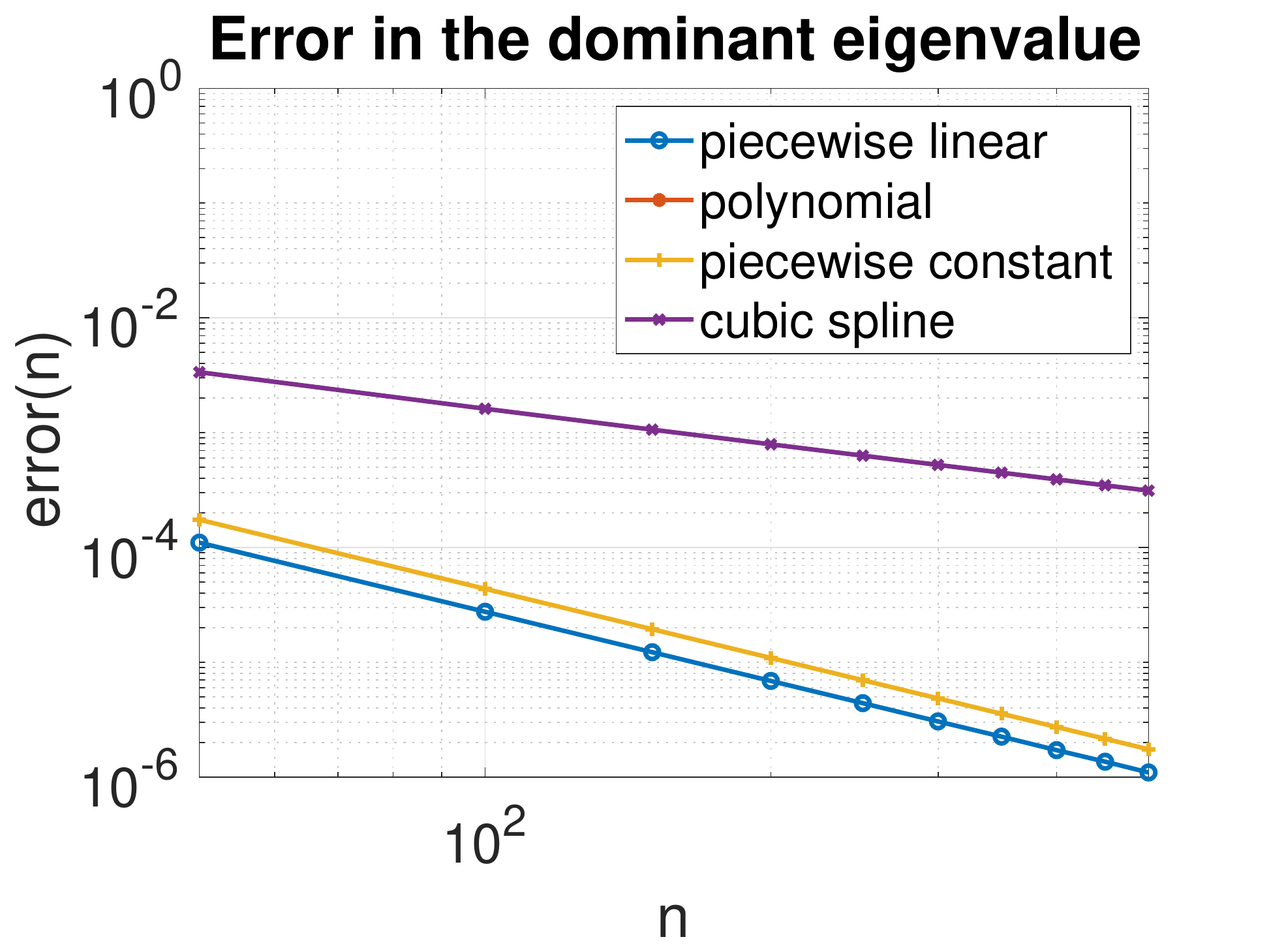}
	\caption{Error in the dominant eigenvalue for piecewise linear collocation $(\circ)$, polynomial collocation ($\bullet$), Bubnov-Galerkin with piecewise constant functions $(+)$ and cubic splines $(\tm)$ up to $n=100$ nodes (left) and up to $n=500$ nodes (right)}
	\label{figcolerror}
\end{figure}

Finally let us test the above methods by means of
\begin{exmp}[Laplace kernel]\label{exlaplace}
	Let $\alpha_1,\alpha_2>0$ denote dispersal rates and a mapping $k:\Omega\tm\Omega\to(0,\infty)$ be continuous. The Fredholm operator
	\begin{align*}
		\eK u
		&:=
		\int_\Omega
		K(\cdot,y)u(y)\d y,&
		K(x,y)&:=
		\begin{pmatrix}
			\tfrac{1}{2\alpha_1}e^{-\tfrac{|x-y|}{\alpha_1}} & -k(x,y)\\
			0 & \tfrac{1}{2\alpha_2}e^{-\tfrac{|x-y|}{\alpha_2}}
		\end{pmatrix}
	\end{align*}
	satisfies Hypothesis $(NL)$ and is compact on $C(\Omega)^2$, as well as on $L^2(\Omega)^2$. The south-east cone $Y_+:=\R_+\tm(-\R_+)$ is spanned by the vectors $e_1:=\binom{1}{0}$, $e_2:=-\binom{0}{1}$ and with $e_i':=e_i$ for $i=1,2$ we consequently obtain
	\begin{align*}
		\iprod{K(x,y)e_1,e_1'}&=\tfrac{1}{2\alpha_1}e^{-\tfrac{|x-y|}{\alpha_1}}>0,&
		\iprod{K(x,y)e_2,e_1'}&=k(x,y)>0,\\
		\iprod{K(x,y)e_1,e_2'}&=0,&
		\iprod{K(x,y)e_2,e_2'}&=\tfrac{1}{2\alpha_2}e^{-\tfrac{|x-y|}{\alpha_2}}>0.
	\end{align*}
	Therefore, the kernel $K(x,y)$ is $Y_+$-positive for all $x,y\in\Omega$. Thus, \tref{thmlinmon} implies that $\eK$ is $C(\Omega)_+^2$-positive, while \tref{thml2} ensures $L^2(\Omega)_+^2$-positivity. With the aid of 
	$$
		\eK_iu(x):=\frac{1}{2\alpha_i}\int_\Omega e^{-\tfrac{|x-y|}{\alpha_i}}u(y)\d y\fall x\in\Omega,\,i=1,2
	$$
	we obtain the spectrum $\sigma(\eK)=\sigma(\eK_1)\cup\sigma(\eK_2)$. 
	\\
	In order to compute the dominant eigenvalue of $\eK$ numerically, we retreat to intervals $\Omega=[-\tfrac{L}{2},\tfrac{L}{2}]$ for some $L>0$. It is shown in \cite[pp.~24--27, Sect.~3.2]{lutscher:19} that the dominant eigenvalue $\lambda^i>0$ of $\eK_i$ and the smallest positive solution $\nu_i$ of the transcendental equation $\tan(\tfrac{L}{2\alpha_i}\nu)=\tfrac{1}{\nu}$ are related via $\lambda^i=\tfrac{1}{1+\nu_i^2}$. We refer to Tab.~\ref{tabev} for exact values and approximate the dominant eigenvalue with positive collocation methods (piecewise linear, polynomial and spline) and positive Bubnov-Galerkin methods (piecewise constant). This is based on a grid \eqref{grid} with $a=-\tfrac{L}{2}$, $b=\tfrac{L}{2}$ and $x_j:=-\tfrac{L}{2}+j\tfrac{L}{n}$, $0\leq j\leq n$. For the sake of discrete projection methods, the remaining integrals are approximated by the summed midpoint rule \eqref{nomid} with the centered nodes $\eta_i:=x_i-\tfrac{L}{2n}$, $1\leq i\leq n$. 
	\\
	For $L=2$ and $\alpha_i$ as in Tab.~\ref{tabev} we moreover compute the dominant eigenvalue of $K^n$ and relate it to the exact eigenvalue $\lambda^1\approx 0.5746552163364324$ of $\eK$. The results of our numerical simulations are illustrated in \fref{figcolerror}. Piecewise linear collocation and the Bubnov-Galerkin method with piecewise constant basis functions illustrate quadratic convergence. This is also true for collocation with polynomial basis functions until beginning with $n\approx 25$ computational instabilities become apparent. The positivity preserving collocation based on cubic splines from \cref{corspline} shows only linear convergence. 
\end{exmp}
\section{Conclusion}
We provided sufficient conditions that positivity properties in a general class of matrix-valued kernels $K$ transfer to Fredholm operators \eqref{nofred} over compact domains $\Omega$. In addition, the persistence of this property under Nystr\"om and projection methods is studied. As a result, when solving numerical problems involving such operators, we recommend to restrict to Nystr\"om methods with positive weights. Among the projection methods, collocation with piecewise linear basis functions yields positivity-preserving semi-discretizations and a combination with the Midpoint Rule \eqref{nomid} or the Trapezoidal Rule \eqref{notrap} leads to a corresponding scheme. Although there are further positivity preserving collocation methods, they have less favorable properties. For Bubnov-Galerkin methods, piecewise constant approximation is suitable and can be combined with e.g.\ the Midpoint Rule \eqref{nomid} in order to arrive at feasible schemes preserving positivity. 
\begin{appendix}
\renewcommand{\theequation}{\Alph{section}.\arabic{equation}}
\section{Cones and positive operators}
\label{secA}
Let $(X,\norm{\cdot})$ denote a real Banach space with dual space $X'$ and the duality pairing $\iprod{x,x'}:=x'(x)$. A nonempty closed and convex subset $X_+\subseteq X$ is called \emph{(order) cone}, if $\R_+X_+\subseteq X_+$ and $X_+\cap(-X_+)=\set{0}$ hold. Equipped with such a cone one arrives at an \emph{ordered Banach space} $X$. Let us assume $X_+\neq\set{0}$ throughout and for elements $x,\bar x\in X$ we write
\begin{align}
	x\leq\bar x\quad&:\Leftrightarrow\quad\bar x-x\in X_+,
	\notag\\
	x<\bar x\quad&:\Leftrightarrow\quad\bar x-x\in X_+\setminus\set{0},
	\label{noorders3}\\
	x\ll \bar x\quad&:\Leftrightarrow\quad\bar x-x\in X_+^\circ; 
	\notag
\end{align}
the latter relation requires $X_+^\circ\neq\emptyset$ and one speaks of a \emph{solid cone} $X_+$. A cone is \emph{total}, if $X=\overline{X_+-X_+}$. Solid cones satisfy $X=X_+-X_+$ and are total. 

By means of the \emph{dual cone} $X_+':=\set{x'\in X':\,0\leq\iprod{x,x'}\text{ for all }x\in X_+}$ we can characterize the elements of $X_+$ and $X_+^\circ$ as follows:
\begin{lem}\label{lemA1}
	\begin{enumerate}
		\item[(a)] $X_+'\neq\{0\}$ and for every $x\in X$ the following holds: 
		\begin{align*}
			x\in X_+&\Leftrightarrow 0\leq\iprod{x,x'}\fall x'\in X_+',\\
			x\in X_+\setminus\set{0}&\Rightarrow 0<\iprod{x,x_0'}\quad\text{for some }x_0'\in X_+'\setminus\set{0}. 
		\end{align*}

		\item[(b)] If $X_+$ is solid, then for every $x\in X$ the following holds:
		\begin{align}
			x\in X_+^\circ&\Leftrightarrow 0<\iprod{x,x'}\fall x'\in X_+'\setminus\set{0},
			\notag\\
			x\in\partial X_+&\Rightarrow 0=\iprod{x,x_0'}\quad\text{for some }x_0'\in X_+'\setminus\set{0}. 
			\label{lemA14}
		\end{align}
	\end{enumerate}
\end{lem}
\begin{pf}
	Due to \cite[p.~222, Prop.~19.3(a--b)]{deimling:85} it remains to verify the direction ``$\Leftarrow$'' in the first equivalence of (b). Let $0<\iprod{x,x'}$ for all $x'\in X_+'\setminus\set{0}$ and we deduce $x\in X_+$ from (a). Using the contraposition of \eqref{lemA14} we derive $x\notin \partial X_+$ which implies that $x\in X_+^\circ$.
	\qed
\end{pf}

A bounded linear mapping $T\in L(X)$ is called 
\begin{itemize}
	\item \emph{positive}, if $T(X_+\setminus\set{0})\subseteq X_+$, 

	\item \emph{strictly positive}, if $T(X_+\setminus\set{0})\subseteq X_+\setminus\set{0}$, 

	\item \emph{strongly positive}, if $T(X_+\setminus\set{0})\subseteq X_+^\circ$.
\end{itemize}
When working with several cones, we write $X_+$-positive etc., in order to indicate a particular cone. We denote $T\in L(X)$ as $X_+$-\emph{injective}, provided
$$
	N(T)\cap X_+=\set{0}
$$
holds. Then $T$ is strictly positive, if and only if it is positive and $X_+$-injective. A strongly positive $T$ yields the inclusion $T X_+^\circ\subseteq X_+^\circ$. 

The positivity properties are preserved under compositions and it holds:
\begin{cor}\label{corA2}
	Let $X_+$ be solid and $T,S\in L(X)$. If $S$ is strictly positive and $T$ is strongly positive, or if $S$ strongly positive and $T$ satisfies $TX_+^\circ\subseteq X_+^\circ$, then $TS$ is strongly positive. 
\end{cor}
\begin{pf}
	This is immediate by definition. 
	\qed
\end{pf}

We conclude with a Krein-Rutman theorem suiting our aims (see also \cite{li:jia:21}):
\begin{thm}[Krein-Rutman]\label{thmkr}
	Let $T\in L(X)$ be compact.
	\begin{enumerate}
		\item[(a)] If $X_+$ is total and $T$ is positive with $r(T)>0$, then the spectral radius $r(T)$ is an eigenvalue of $T$ with positive eigenvector $x^\ast$ (cf.\ \cite[p.~226, Thm.~19.2]{deimling:85} or \cite[p.~290, Prop.~7.26]{zeidler:93}).

		\item[(b)] If $X_+$ is solid and $T$ is strongly positive, then $T$ has exactly one eigenvector with $0<x^\ast$ and $\norm{x^\ast}=1$; the corresponding eigenvalue is $r(T)$ and $0\ll x^\ast$ (cf.\ \cite[p.~290, Thm.~7.C]{zeidler:93}).
	\end{enumerate}
	One calls $x^\ast$ a \emph{dominant eigenvector} (\emph{eigenfunction} on function spaces $X$). 
\end{thm}
\section{Tridiagonal matrices}
Let $n\in\N\setminus\set{1}$ and $a_1,\ldots,a_n$, $b_1,\ldots,b_{n-1}$, $c_1,\ldots,c_{n-1}$ be reals with $b_j\neq 0$, $1\leq j<n$, such that the tridiagonal matrix
$$
	T:=
	\begin{pmatrix}
		a_1 & b_1 & & &\\
		c_1 & a_2 & b_2 & &\\
		& \ddots & \ddots & \ddots & \\
		&& c_{n-2} & a_{n-1} & b_{n-1} \\
		&&& c_{n-1} & a_n
	\end{pmatrix}
$$
is nonsingular. If we recursively introduce the finite sequences
\begin{align*}
	d_n&:=a_n,&
	d_{t-1}&:=a_{t-1}-\tfrac{b_{t-1}c_{t-1}}{d_t}\fall t=n,\ldots,2,\\
	\delta_1&:=a_1,&
	\delta_{t+1}&:=a_{t+1}-\tfrac{b_tc_t}{\delta_t}\fall t=1,\ldots,n-1,
\end{align*}
then the following holds (with the convention that empty products are $1$): 
\begin{lem}[{\cite[Thm.~2.1]{fonseca:petronilho:01}}]\label{leminv}
	The entries of the inverse $T^{-1}$ are given by
	$$
		(T^{-1})_{ij}
		=
		(-1)^{i+j}
		\begin{cases}
			b_i\cdots b_{j-1}\tfrac{d_{j+1}\cdots d_n}{\delta_i\cdots\delta_n},&i\leq j,\\
			c_j\cdots c_{i-1}\tfrac{d_{i+1}\cdots d_n}{\delta_j\cdots\delta_n},&j<i.
		\end{cases}
	$$
\end{lem}
\section{Quadrature rules}
\label{secC}
Below we list the quadrature rules \eqref{quad} used in our numerical simulations. We abbreviate $x_j:=a+jh$ with $h=\tfrac{b-a}{n}$ and refer to \cite[pp.~361ff, Chap.~15]{engeln:uhlig:96} for the following facts (here $\xi\in[a,b]$ refers to an intermediate point).
\begin{itemize}
	\item Summed midpoint rule: The constant weights $w_j=h$ supplemented by nodes $\eta_j=a+(j+\tfrac{1}{2})h$ and $q_n:=n-1$ lead to
	\begin{equation}
		\int_a^bu
		=
		h\sum_{j=0}^{n-1}u(x_j+\tfrac{h}{2})+\tfrac{b-a}{24}h^2u''(\xi)
		\label{nomid}
	\end{equation}

	\item Summed trapezoidal rule: The weights $w_0=w_n=\tfrac{h}{2}$, $w_j=h$ for $1\leq j<n$, the nodes $\eta_j=x_j$ and $q_n=n$ yield
	\begin{equation}
		\int_a^bu
		=
		\tfrac{h}{2}\sum_{j=0}^{n-1}\intoo{u(x_j)+u(x_{j+1})}-\tfrac{b-a}{12}h^2u''(\xi)
		\label{notrap}
	\end{equation}

	\item Summed Milne's rule: The weights $w_{3j}=w_{3j+2}=\tfrac{2}{3}h$, $w_{3j+2}=-\tfrac{1}{3}h$, the nodes $\eta_{3j}=\tfrac{3x_j+x_{j+1}}{4}$, $\eta_{3j+1}=\tfrac{x_j+x_{j+1}}{2}$, $\eta_{3j+2}=\tfrac{x_j+3x_{j+1}}{4}$ and $q_n=3n$ give
\begin{align}
	\int_a^bu
	=&
	\tfrac{h}{3}\sum_{j=0}^{n-1}\intoo{2u(\tfrac{3x_j+x_{j+1}}{4})-u(\tfrac{x_j+x_{j+1}}{2})+2u(\tfrac{x_j+3x_{j+1}}{4})}
	\notag\\
	&+\tfrac{b-a}{23040}h^4u^{(4)}(\xi)
	\label{nomilne}
\end{align}

	\item Summed $6$th order Gau{\ss}: The weights $w_{3j}=w_{3j+2}=\tfrac{5}{8}h$, $w_{3j+1}=\tfrac{4}{9}h$, the nodes 
	$\eta_{3j}=x_j+\bigl(1-\sqrt{\tfrac{3}{5}}\bigr)\tfrac{h}{2}$, 
	$\eta_{3j+1}=x_j+\tfrac{h}{2}$, 
	$\eta_{3j+2}=x_j+\bigl(1+\sqrt{\tfrac{3}{5}}\bigr)\tfrac{h}{2}$ and $q_n=3n$ imply the quadrature rule
	\begin{align}
		\int_a^bu
		=&
		\tfrac{h}{18}\sum_{j=0}^{n-1}\intoo{
		5u(x_j+(1-\sqrt{\tfrac{3}{5}})\tfrac{h}{2})+
		8u(x_j+\tfrac{h}{2})+
		5u(x_j+(1+\sqrt{\tfrac{3}{5}})\tfrac{h}{2})}
		\notag\\
		&
		+\tfrac{b-a}{31500}h^6u^{(6)}(\xi).
		\label{nogauss}
	\end{align}
\end{itemize}
\end{appendix}
%
%
%
%
%
%
%

%
%
%
\end{document}